\documentclass[11pt]{article}

\usepackage{amsmath}
\usepackage{amssymb}
\usepackage{amsthm}
\usepackage{latexsym}
\usepackage{color}
\usepackage{graphicx}
\usepackage{color}
\usepackage[T1]{fontenc}
\usepackage[english]{babel}
\usepackage{geometry}
\DeclareSymbolFont{calletters}{OMS}{cmsy}{m}{n}
\DeclareSymbolFontAlphabet{\mathcal}{calletters}

\def\be{\begin{eqnarray}}
\def\ee{\end{eqnarray}}

\def\b*{\begin{eqnarray*}}
\def\e*{\end{eqnarray*}}

\newtheorem{Theorem}{Theorem}[part]
\newtheorem{Definition}{Definition}[part]
\newtheorem{Proposition}{Proposition}[part]

\newtheorem{Assumption}{Assumption}[part]
\newtheorem{Lemma}{Lemma}[part]

\newtheorem{Remark}{Remark}[part]

\makeatletter \@addtoreset{equation}{section}

\@addtoreset{Definition}{section}

\@addtoreset{Theorem}{section}

\@addtoreset{Proposition}{section}

\@addtoreset{Property}{section}

\@addtoreset{Assumption}{section}

\@addtoreset{Corollary}{section}

\@addtoreset{Lemma}{section}

\@addtoreset{Remark}{section}

\@addtoreset{Example}{section}

\newcommand{\No}[1]{\left\|#1\right\|}     
\newcommand{\abs}[1]{\left|#1\right|}     

\def \F{\mathbb{F}}
\def \H{\mathbb{H}}

\def \P{\mathbb{P}}

\def\Fc{{\cal F}}

\def\Tr#1{{\rm Tr}\left[#1\right]}

\def \Sup{\displaystyle\sup}

\def\einf{{\rm ess \, inf}}
\def\esup{{\rm ess \, sup}}
\def\trace{{\rm Tr}}

\def\={\;=\;}
\def\.{\;.}

\def\eps{\varepsilon}

\def\reff#1{{\rm(\ref{#1})}}

\def\1{{\bf 1}}
\def \ep{\hbox{ }\hfill{ ${\cal t}$~\hspace{-5.1mm}~${\cal u}$   } }

\def \proof{{\noindent \bf Proof. }}
\def \ep{\hbox{ }\hfill$\Box$}

 \def\normeL2#1{\left\|{#1}\right\|_{L^2}}

 \title{Second Order Backward Stochastic Differential Equations with Quadratic Growth\footnote{Research supported by the Chair {\it Financial Risks} of the {\it Risk Foundation} sponsored by Soci\'et\'e G\'en\'erale, the Chair {\it Derivatives of the Future} sponsored by the {F\'ed\'eration Bancaire Fran\c{c}aise}, and the Chair {\it Finance and Sustainable Development} sponsored by EDF and Calyon.}
}
\author{Dylan {\sc Possama\"{i}}\footnote{CMAP, Ecole Polytechnique, Paris, dylan.possamai@polytechnique.edu.}
      \and Chao {\sc Zhou}\footnote{CMAP, Ecole Polytechnique, Paris, chao.zhou@polytechnique.edu.} \footnote{The authors are grateful to Anis Matoussi, Guillaume Royer, Xiaolu Tan and Nizar Touzi for their help and precious remarks. The authors would also like to thank two anonymous referees and an associate editor, whose advices greatly helped us to improve a previous version of this paper.}
}
 \date{\today}

\begin{document}

 \maketitle

\vspace{3mm}

 \begin{abstract}
We extend the wellposedness results for second order backward stochastic differential equations introduced by Soner, Touzi and Zhang \cite{stz} to the case of a bounded terminal condition and a generator with quadratic growth in the $z$ variable. More precisely, we obtain uniqueness through a representation of the solution inspired by stochastic control theory, and we obtain two existence results using two different methods. In particular, we obtain the existence of the simplest purely quadratic 2BSDEs through the classical exponential change, which allows us to introduce a quasi-sure version of the entropic risk measure. As an application, we also study robust risk-sensitive control problems. Finally, we prove a Feynman-Kac formula and a probabilistic representation for fully non-linear PDEs in this setting.
\vspace{10mm}

\noindent{\bf Key words:} Second order backward stochastic differential equation, quadratic growth, BMO martingales, r.c.p.d., Feyman-Kac, fully non-linear PDEs, quasi-sure.
\vspace{5mm}

\noindent{\bf AMS 2000 subject classifications:} 60H10, 60H30
\end{abstract}
\newpage

\section{Introduction}

Backward stochastic differential equations (BSDEs for short) appeared for the first time in Bismut \cite{bis} in the linear case. However, they only became a popular field of research after the seminal paper of Pardoux and Peng \cite{pardpeng}, mainly because of the very large scope of their domain of applications, ranging from stochastic control to mathematical finance.

\vspace{0.2em}
On a filtered probability space $(\Omega,\mathcal F,\left\{\mathcal F_t\right\}_{0\leq t\leq T},\mathbb P)$ generated by an $\mathbb R^d$-valued Brownian motion $B$, solving a BSDE amounts to finding a pair of progressively measurable processes $(Y,Z)$ such that
\begin{equation*}
Y_t=\xi -\int_t^Tf_s(Y_s,Z_s)ds-\int_t^TZ_sdB_s,\text{ } t\in [0,T], \text{ }\mathbb P-a.s.
\end{equation*}
where $f$ (called the generator) is a progressively measurable function and $\xi$ is an $\mathcal F_T$-measurable random variable.

\vspace{0.2em}
Pardoux and Peng \cite{pardpeng} proved existence and uniqueness of the above BSDE provided that the function $f$ is uniformly Lipschitz in $y$ and $z$ and that $\xi$ and $f_s(0,0)$ are square integrable. Then, in \cite{pardpeng2}, they proved that if the randomness in $f$ and $\xi$ is induced by the current value of a state process defined by a forward stochastic differential equation, then the solution to the BSDE could be linked to the solution of a semi-linear PDE by means of a generalized Feynman-Kac formula. This link then opened the way to probabilistic numerical methods for solving semi-linear PDEs, which particularly well suited for highly dimensional problems (see Bouchard and Touzi \cite{bt} among many other).

\vspace{0.2em}
Nonetheless, the class of fully non-linear PDEs remained inaccessible when considering only BSDEs, until a recent work of Cheredito, Soner, Touzi and Victoir \cite{cstv}. They introduced a notion of second order BSDEs ($2$BSDEs), which were then proved to be naturally linked to fully non-linear PDEs. However, only a uniqueness result in the Markovian case was proved, and no existence result (apart from trivial ones) were available. Following this work, Soner, Touzi and Zhang \cite{stz} gave a new definition of 2BSDEs, and thus provided a complete theory of existence and uniqueness under uniform Lipschitz conditions similar to those of Pardoux and Peng. Their key idea was to reinforce the condition that the $2$BSDE must hold $\mathbb P-a.s.$ for every probability measure $\mathbb P$ in a non-dominated class of mutually singular measures (see Section \ref{section.0quad} for precise definitions). Let us describe briefly the intuition behind their definition.

\vspace{0.2em}
Suppose that we want to study the following fully non-linear PDE
\begin{equation}
- \frac{\partial u}{\partial t}-h\left(t,x,u(t,x),Du(t,x),D^2u(t,x)\right)=0,\qquad u(T,x)=g(x).
\label{eq:1}
\end{equation}

If the function $\gamma\mapsto h(t,x,r,p,\gamma)$ is assumed to be convex, then it is equal to its double Fenchel-Legendre transform, and if we denote its Fenchel-Legendre transform by $f$, we have
\begin{equation}
h(t,x,r,p,\gamma)=\underset{a\geq 0}{\sup}\left\{\frac12a\gamma-f(t,x,r,p,a)\right\}
\label{eq:2}
\end{equation}

Then, from \reff{eq:2}, we expect, at least formally, that the solution $u$ of \reff{eq:1} is going to verify
$$u(t,x)=\underset{a\geq 0}{\sup}\ u^a(t,x),$$
where $u^a$ is defined as the solution of the following semi-linear PDE
\begin{equation}
- \frac{\partial u^a}{\partial t}-\frac12aD^2u^a(t,x)+f\left(t,x,u^a(t,x),Du^a(t,x),a\right)=0,\qquad u^a(T,x)=g(x).
\label{eq:3}
\end{equation}

Since $u^a$ is linked to a classical BSDE, the 2BSDE associated to $u$ should correspond (in some sense) to the supremum of the family of BSDEs indexed by $a$. Furthermore, changing the process $a$ can be achieved by changing the probability measure under which the BSDE is written. However, this amounts to changing the quadratic variation of the martingale driving the BSDE, and therefore leads to a family of mutually singular probability measures. In these respects, the 2BSDE theory shares deep links with the theory of quasi-sure stochastic analysis of Denis and Martini \cite{denis} and the theory of $G$-expectation of Peng \cite{peng}.

\vspace{0.2em}
Following the breakthrough of \cite{stz}, Possamai \cite{pos} extended the existence and uniqueness result for 2BSDEs to the case of a generator having linear growth and satisfying a monotonicity condition. Motivated by a robust utility maximization problem under volatility uncertainty (see the accompanying paper \cite{mpz}), our aim here is to go beyond the results of \cite{pos} to prove an existence and uniqueness result for $2$BSDEs whose generator has quadratic growth in $z$.

\vspace{0.2em}
The question of existence and uniqueness of solutions to these quadratic equations in the classical case was first examined by Kobylanski \cite{kob}, who proved existence and uniqueness of a solution by means of approximation techniques borrowed from the PDE literature, when the generator is continuous and has quadratic growth in $z$ and the terminal condition $\xi$ is bounded. Then, Tevzadze \cite{tev} has given a direct proof for the existence and uniqueness of a bounded solution in the Lipschitz-quadratic case, proving the convergence of the usual Picard iteration. Following those works, Briand and Hu \cite{bh} have extended the existence result to unbounded terminal condition with exponential moments and proved uniqueness for a convex coefficient \cite{bh2}. Finally, Barrieu and El Karoui \cite{elkarbar} recently adopted a completely different approach, embracing a forward point of view to prove existence under conditions similar to those of Briand and Hu.

\vspace{0.2em}
With quadratic growth generators, the proof of wellposedness of 2BSDEs becomes more technical. In this paper, we propose two very different methods to prove it. Our main contributions are the following. First, as for classical quadratic BSDEs, we find a link between quadratic 2BSDEs and a generalized BMO space. Moreover, in contrast with the classical framework, for which the only property of BMO martingales used is the fact that their Dol\'eans-Dade exponential is a uniformly integrable martingale, in the 2BSDE framework, we use extensively some fine properties of BMO martingales. This seems to be directly linked to the presence of the non-decreasing processes in the definition of 2BSDEs. With these estimates in hand, we manage to recover a uniqueness result using a representation of the solution similar to the ones obtained previously in the literature. Concerning our proofs of existence, one is inspired by the pathwise construction initiated in \cite{stz}. In this regard, the proof is very similar. However, our contribution lies the fact that on the one hand we point out some sufficient properties which must be satisfied in order to apply this type of proof, and on the other hand, we prove a new technical result (see Proposition \ref{prop.tech}) which allows to obtain the existence without using the Picard argument anymore as in \cite{stz}. Next, we investigate the approximation approach for the existence of a solution as for the classical BSDEs, and discuss in detail why in general this approach fails. However, we still manage to find some cases where it works and thus provide a second existence result. Finally, keeping in mind the possible applications of our theoretical results (see also \cite{mpz}), we define and solve some robust risk-sensitive control problems using 2BSDEs with quadratic growth, and we show that the links with fully non-linear PDEs proved in \cite{stz} still hold in our context.

\vspace{0.2em}
The rest of the paper is organized as follows. First, we recall some notations in Section \ref{section.0quad} and prove a uniqueness result in Section \ref{uni} by means of a priori estimates and a representation of the solution inspired by the stochastic control theory. Section \ref{sec.pure} is devoted to the study of purely quadratic 2BSDEs and robust risk-sensitive control problems. Next, in Section \ref{sec.approx}, we consider and discuss approximation techniques for the problem of existence of a solution. Then, in Section \ref{section.1quad}, we use a completely different method introduced by Soner, Touzi and Zhang \cite{stz} to construct the solution to the quadratic $2$BSDE path by path. Finally, in Section \ref{section.3quad}, we extend the results of Soner, Touzi and Zhang \cite{stz} on the connections between fully non-linear PDEs and $2$BSDEs to the quadratic case.


\section{Preliminaries} \label{section.0quad}

Let $\Omega:=\left\{\omega\in C([0,T],\mathbb R^d):\omega_0=0\right\}$ be the canonical space equipped with the uniform norm $\No{\omega}_{\infty}:=\sup_{0\leq t\leq T}|\omega_t|$, $B$ the canonical process, $\mathbb P_0$ the Wiener measure, $\mathbb F:=\left\{\mathcal F_t\right\}_{0\leq t\leq T}$ the filtration generated by $B$, and $\mathbb F^+:=\left\{\mathcal F_t^+\right\}_{0\leq t\leq T}$ the right limit of $\mathbb F$.

\subsection{A first set of probability measures}\label{subsec.proba}

Our aim here is to give a correct mathematical basis to the intuitions we provided in the Introduction. We first recall that by the results of Karandikar \cite{kar}, we can give pathwise definitions of the quadratic variation $\left<B\right>_t$ and its density $\widehat a_t$. Let $\overline{\mathcal P}_W$ denote the set of all local martingale measures $\mathbb P$ (i.e. the probability measures $\mathbb P$ under which $B$ is a local martingale) such that $\left<B\right>_t \text{ is absolutely continuous in $t$ and $\widehat a$ takes values in $\mathbb S_d^{>0}$, }\mathbb P-a.s.$, where $\mathbb S_d^{>0}$ denotes the space of all $d\times d$ real valued symmetric positive definite matrices. As in \cite{stz}, we concentrate on the subclass $\overline{\mathcal P}_S\subset\overline{\mathcal P}_W$ consisting of all probability measures
\begin{equation}\label{proba.pquad}
\mathbb P^\alpha:=\mathbb P_0\circ (X^\alpha)^{-1} \text{ where }X_t^\alpha:=\int_0^t\alpha_s^{1/2}dB_s,\text{ }t\in [0,T],\text{ }\mathbb P_0-a.s.
\end{equation}
for some $\mathbb F$-progressively measurable process $\alpha$ satisfying $\int_0^T\abs{\alpha_s}ds<+\infty$. We recall from \cite{stz2} that every $\mathbb P\in \overline{\mathcal P}_S$ satisfies the Blumenthal zero-one law and the martingale representation property. Notice that the set $\overline{\mathcal P}_S$ is bigger than the set $\widetilde{\mathcal P}_S$ introduced in \cite{pos}, which is defined,for some fixed matrices $\underline a$ and $\bar a$ in $\mathbb S_d^{>0}$, by
\begin{equation}\label{ps}
\widetilde{\mathcal P}_S:=\left\{\mathbb P^\alpha\in\overline{\mathcal P}_S,\ \underline a\leq\alpha\leq\bar a, \ \mathbb P_0-a.s.\right\},
\end{equation}

\subsection{The Generator and the final set $\mathcal P_H$}
Before defining the spaces under which we will be working or defining the 2BSDE itself, we first need to restrict one more time our set of probability measures, using explicitly the generator of the 2BSDE. Following the PDE intuition, let us first consider a map $H_t(\omega,y,z,\gamma):[0,T]\times\Omega\times\mathbb{R}\times\mathbb{R}^d\times D_H\rightarrow \mathbb{R}$, where $D_H \subset \mathbb{R}^{d\times d}$ is a given subset containing $0$. As expected, we define its Fenchel-Legendre conjugate w.r.t.$\gamma$ by
\begin{align*}
&F_t(\omega,y,z,a):=\underset{\gamma \in D_H}{\Sup}\left\{\frac12\trace(a\gamma)-H_t(\omega,y,z,\gamma)\right\} \text{ for } a \in \mathbb S_d^{>0}\\[0.3em]
&\widehat{F}_t(y,z):=F_t(y,z,\widehat{a}_t) \text{ and } \widehat{F}_t^0:=\widehat{F}_t(0,0).
\end{align*}

We denote by $D_{F_t(y,z)}$ the domain of $F$ in $a$ for a fixed $(t,\omega,y,z)$, and as in \cite{stz} we restrict the probability measures in $\mathcal{P}_H\subset \overline{\mathcal{P}}_S$

\begin{Definition}\label{defquad}
$\mathcal{P}_H$ consists of all $\mathbb P \in \overline{\mathcal{P}}_S$ such that
$$\underline{a}_\mathbb P\leq \widehat{a}\leq \bar{a}_\mathbb P, \text{ } dt\times d\mathbb P-a.s. \text{ for some } \underline{a}_\mathbb P, \bar{a}_\mathbb P \in \mathbb S_d^{>0}, \text{and $\widehat a_t\in D_{F_t(y,z)}$. }$$
\end{Definition}

\begin{Remark}
The restriction to the set $\mathcal P_H$ obeys two imperatives. First, since $\widehat F$ is destined to be the generator of our 2BSDE, we obviously need to restrict ourselves to probability measures such that $\widehat a_t\in D_{F_t(y,z)}$. Moreover, we also restrict the measures considered to the ones such that the density of the quadratic variation of $B$ is bounded to ensure that $B$ is actually a true martingale under each of those probability measures. This will be important to obtain a priori estimates.
\end{Remark}

\subsection{Spaces of interest}

We now provide the definition of the spaces which will be used throughout the paper.
\subsubsection{Quasi-sure spaces}\label{sec.space}
$\mathbb L^{\infty}_H$ denotes the space of all $\mathcal F_T$-measurable scalar r.v. $\xi$ with
$$\No{\xi}_{\mathbb L^{\infty}_H}:=\underset{\mathbb{P} \in \mathcal{P}_H}{\sup}\No{\xi}_{L^{\infty}(\mathbb P)}<+\infty.$$

$\mathbb H^{p}_H$ denotes the space of all $\mathbb F^+$-progressively measurable $\mathbb R^d$-valued processes $Z$ with
$$\No{Z}_{\mathbb H^{p}_H}^p:=\underset{\mathbb{P} \in \mathcal{P}_H}{\sup}\mathbb E^{\mathbb P}\left[\left(\int_0^T|\widehat a_t^{1/2}Z_t|^2dt\right)^{\frac p2}\right]<+\infty.$$

$\mathbb D^{\infty}_H$ denotes the space of all $\mathbb F^+$-progressively measurable $\mathbb R$-valued processes $Y$ with
$$\mathcal P_H-q.s. \text{ c\`adl\`ag paths, and }\No{Y}_{\mathbb D^{\infty}_H}:=\underset{0\leq t\leq T}{\sup}\No{Y_t}_{\mathbb L^{\infty}_H}<+\infty.$$

Finally, we denote by $\mbox{UC}_b(\Omega)$ the collection of all bounded and uniformly continuous maps $\xi:\Omega\rightarrow \mathbb R$ with respect to the $\No{\cdot}_{\infty}$-norm, and we let
$$\mathcal L^{\infty}_H:=\text{the closure of $\mbox{UC}_b(\Omega)$ under the norm $\No{\cdot}_{\mathbb L^{\infty}_H}$}.$$

\begin{Remark}
We emphasize that all the above norms and spaces are the natural generalizations of the standard ones when working under a single probability measure.
\end{Remark}

\subsubsection{The space $\mathbb B\rm{MO}(\mathcal P_H)$ and important properties}

It is a well known fact that the $Z$ component of the solution of a quadratic BSDE with a bounded terminal condition belongs to the so-called BMO space. Since this link will be extended and used intensively throughout the paper, we will recall some results and definitions for the BMO space, and then extend them to our quasi-sure framework. We first recall (with a slight abuse of notation) the definition of the BMO space for a given probability measure $\mathbb P$.

\begin{Definition}
$\rm{BMO}(\mathbb P)$ denotes the space of all $\mathbb F^+$-progressively measurable $\mathbb R^d$-valued processes $Z$ with
$$\No{Z}_{\rm{BMO}(\mathbb P)}:=\underset{\tau\in\mathcal T_0^T}{\sup}\No{\mathbb E_\tau^\mathbb P\left[\int_\tau^T|\widehat a_t^{1/2}Z_t|^2dt\right]}_{\mathbb L^{\infty}(\mathbb P)}<+\infty,$$
where $\mathcal T_0^T$ is the set of $\mathcal F_t$ stopping times taking their values in $[0,T]$.
\end{Definition}

We also recall the so called energy inequalities (see \cite{kaz} and the references therein). Let $Z\in\rm{BMO}(\mathbb P)$ and $p\geq 1$. Then we have
\begin{align}\label{energy}
\mathbb E^\mathbb P\left[\left(\int_0^T\abs{\widehat a_s^{1/2}Z_s}^2ds\right)^p\right]\leq 2p!\left(4\No{Z}_{\rm{BMO}(\mathbb P)}^2\right)^p.
\end{align}

The extension to a quasi-sure framework is then naturally given by the following space.

\vspace{0.2em}
$\mathbb B\rm{MO}(\mathcal P_H)$ denotes the space of all $\mathbb F^+$-progressively measurable $\mathbb R^d$-valued processes $Z$ with
$$\No{Z}_{\mathbb B\rm{MO}(\mathcal P_H)}:=\underset{\mathbb{P} \in \mathcal{P}_H}{\sup}\No{Z}_{\rm{BMO}(\mathbb P)}<+\infty.$$

The main interest of the BMO spaces is that if a process $Z$ belongs to it, then the stochastic integral $\int_0^.Z_sdB_s$ is a uniformly integrable martingale, which in turn allows us to use it for changing the probability measure considered via Girsanov's Theorem. The two following results give more detailed results in terms of $L^r$ integrability of the corresponding Dol\'eans-Dade exponential.

\begin{Lemma}\label{lemma.bmoholder}
Let $Z\in\mathbb B\rm{MO}(\mathcal P_H)$. Then there exists $r>1$, such that
$$\underset{\mathbb P\in\mathcal P_H}{\Sup}\mathbb E^\mathbb P\left[\left(\mathcal E\left(\int_0^.Z_sdB_s\right)\right)^r\right]<+\infty.$$
\end{Lemma}

\vspace{0.2em}
\proof
By Theorem $3.1$ in \cite{kaz}, we know that if $\No{Z}_{\rm{BMO}(\mathbb P)}\leq \Phi(r)$ for some one-to-one function $\Phi$ from $(1,+\infty)$ to $\mathbb R^*_+$, then $\mathcal E\left(\int_0^.Z_sdB_s\right)$ is in $L^r(\mathbb P)$. Here, since $Z\in\mathbb B\rm{MO}(\mathcal P_H)$, the same $r$ can be used for all the probability measures.
\ep

\vspace{0.2em}
\begin{Lemma}\label{lemma.bmo2}
Let $Z\in\mathbb B\rm{MO}(\mathcal P_H)$. Then there exists $r>1$, such that for all $t\in[0,T]$
 $$\underset{\mathbb P\in\mathcal P_H}{\sup}\mathbb E_t^\mathbb P\left[\left(\frac{\mathcal E\left(\int_0^tZ_sdB_s\right)}{\mathcal E\left(\int_0^TZ_sdB_s\right)}\right)^{\frac{1}{r-1}}\right]<+\infty.$$
\end{Lemma}

\vspace{0.2em}
\proof
This is a direct application of Theorem $2.4$ in \cite{kaz} for all $\mathbb P\in\mathcal P_H$.
\ep

\vspace{0.2em}
We emphasize that the two previous Lemmas are absolutely crucial to our proof of uniqueness and existence. Besides, they also play a major role in our accompanying paper \cite{mpz}.

\subsection{The definition of the 2BSDE}
Everything is now ready to define the solution of a 2BSDE. We shall consider the following
\begin{equation}
Y_t=\xi -\int_t^T\widehat{F}_s(Y_s,Z_s)ds -\int_t^TZ_sdB_s + K_T-K_t, \text{ } 0\leq t\leq T, \text{  } \mathcal P_H-q.s.
\label{2bsdequad}
\end{equation}

\begin{Definition}
We say $(Y,Z)\in \mathbb D^{\infty}_H\times\mathbb H^{2}_H$ is a solution to $2$BSDE \reff{2bsdequad} if
\begin{itemize}
\item[$\bullet$] $Y_T=\xi$, $\mathcal P_H-q.s.$
\item[$\bullet$] For all $\mathbb P \in \mathcal P_H$, the process $K^{\mathbb P}$ defined below has non-decreasing paths $\mathbb P-a.s.$
\begin{equation}
K_t^{\mathbb P}:=Y_0-Y_t + \int_0^t\widehat{F}_s(Y_s,Z_s)ds+\int_0^tZ_sdB_s, \text{ } 0\leq t\leq T, \text{  } \mathbb P-a.s.
\label{2bsde.kquad}
\end{equation}

\item[$\bullet$] The family $\left\{K^{\mathbb P}, \mathbb P \in \mathcal P_H\right\}$ satisfies the minimum condition
\begin{equation}
K_t^{\mathbb P}=\underset{ \mathbb{P}^{'} \in \mathcal{P}_H(t^+,\mathbb{P}) }{ \einf^{\mathbb P} }\mathbb{E}_t^{\mathbb P^{'}}\left[K_T^{\mathbb{P}^{'}}\right], \text{ } 0\leq t\leq T, \text{  } \mathbb P-a.s., \text{ } \forall \mathbb P \in \mathcal P_H,
\label{2bsde.minquad}
\end{equation}
where $\mathcal{P}_H(t^+,\mathbb{P}):=\left\{\mathbb P^{'}\in\mathcal P_H,\ \mathbb P^{'}=\mathbb P\text{ on }\mathcal F_{t^+}\right\}$.
\end{itemize}
Moreover if the family $\left\{K^{\mathbb P}, \mathbb P \in \mathcal P_H\right\}$ can be aggregated into a universal process $K$, we call $(Y,Z,K)$ a solution of $2$BSDE \reff{2bsdequad}.
\end{Definition}

\begin{Remark}\label{rem.int}
Let us comment on this definition. As already explained, the PDE intuition leads us to think that the solution of a 2BSDE should be a supremum of solution of standard BSDEs. Therefore for each $\mathbb P$, the role of the non-decreasing process $K^\mathbb P$ is in some sense to "push" the process $Y$ to remain above the solution of the BSDE with terminal condition $\xi$ and generator $\widehat F$ under $\mathbb P$. In this regard, 2BSDEs share some similarities with reflected BSDEs. Pursuing this analogy, the minimum condition \reff{2bsde.minquad} tells us that the processes $K^\mathbb P$ act in a "minimal" way (exactly as implied by the Skorohod condition for reflected BSDEs), and we will see in the next Section that it implies uniqueness of the solution. Besides, if the set $\mathcal P_H$ was reduced to a singleton
$\{\mathbb P\}$, then \reff{2bsde.minquad} would imply that $K^\mathbb P$ is a martingale and a non-decreasing process and is therefore null. Thus we recover the standard BSDE theory. Next, we would like to emphasize that in the language of G-expectation of Peng \cite{peng}, (\ref{2bsde.minquad}) is equivalent, at least if the family can be aggregated into a process $K$, to saying that $-K$ is a G-martingale. This link has already been observed in \cite{stz3} where the authors proved the G-martingale representation property, which formally corresponds to a 2BSDE with a generator equal to $0$.
\end{Remark}

\begin{Remark}
Concerning the aggregation of the family $\left(K^\mathbb P\right)_{\mathbb P\in\mathcal P_H}$, it was shown in a similar context in Matoussi, Possama\"i and Zhou \cite{mpz2} (see Theorem $4.1$), using recent results of Nutz \cite{nutz}, that it could always be aggregated into a universal process $K$ under additional assumptions (related to axiomatic set theory). We refer the reader to the statement of Theorem \ref{mainquad} for more details.
\end{Remark}

\subsection{Assumptions}\label{assump}
We finish this Section by giving our main assumptions on the generator $\widehat F$ of the 2BSDE.

\begin{Assumption} \label{assump.hquad}
\begin{itemize}
\item[\rm{(i)}] $\mathcal P_H$ is not empty and $D_{F_t(y,z)}=D_{F_t}$ is independent of $(\omega,y,z)$.
\item[\rm{(ii)}] For fixed $(y,z,\gamma)$, $F$ is $\mathbb{F}$-progressively measurable.
\item[\rm{(iii)}] $F$ is uniformly continuous in $\omega$ for the $||\cdot||_\infty$ norm.
\item[\rm{(iv)}] $F$ is continuous in $z$ and there exist $(\alpha, \beta,\gamma)\in \mathbb R_+\times\mathbb R_+\times \mathbb R^*_+$ such that
\begin{equation*}
\abs{ F_t(\omega,y,z,a)}\leq \alpha+\beta\abs{y}+\frac\gamma2\abs{ a^{1/2}z}^2,\text{ for all }(t,y,z,\omega,a).
\end{equation*}
\item[\rm{(v)}] There exist $\mu>0$ and $\phi\in\mathbb B\rm{MO}(\mathcal P_H)$ such that for all $(t,y,z,z',\omega,a),$
\begin{equation*}
\abs{ F_t(\omega,y,z,a)- F_t(\omega,y,z',a)-\phi_t. a^{1/2}(z-z')}\leq \mu a^{1/2}\abs{z-z'}\left(\abs{a^{1/2}z}+\abs{ a^{1/2}z'}\right).
\end{equation*}
\item[\rm{(vi)}] We have the following uniform Lipschitz-type property in $y$
\begin{equation*}
\abs{F_t(\omega,y,z,a)- F_t(\omega,y',z,a)}\leq C\abs{y-y'}, \text{ for all }(y,y',z,t,\omega,a).
\end{equation*}
\end{itemize}
\end{Assumption}

\begin{Remark}\label{remark.debutquad}
Let us comment on the above assumptions. Assumptions \ref{assump.hquad} $\rm{(i)}$ and $\rm{(iii)}$ are taken from \cite{stz} and are needed to deal with the technicalities induced by the quasi-sure framework. Moreover, since $F$ is defined as the conjugate of $H$, it is easy to verify that a sufficient condition for the domain of $F$ to not depend on $(\omega,y,z)$ is that the function $H$ itself verifies the (local) Lipschitz conditions (v) and (vi), as well as being uniformly continuous in $\omega$.   Assumptions \ref{assump.hquad} $\rm{(ii)}$, $\rm{(iv)}$,$\rm{(v)}$ and $\rm{(vi)}$ are quite standard in the classical BSDE literature. By Assumption \ref{assump.hquad}$\rm{(iv)}$, we have that $\widehat F^0_t$ is actually bounded, so we do not need to assume the strong integrability condition introduced in \cite{stz}.
\end{Remark}

\section{A priori estimates and uniqueness of the solution}\label{uni}

Before proving some a priori estimates for the solution of the 2BSDE \reff{2bsdequad}, we will first prove rigorously the intuition given in the Introduction saying that the solution of the 2BSDE should be, in some sense, a supremum of solution of standard BSDEs.
Hence, for any $\mathbb{P}\in\mathcal{P}_H$, $\mathbb{F}$-stopping time $\tau$, and $\mathcal{F}_\tau$-measurable random variable $\xi \in \mathbb{L}^\infty(\mathbb P)$, we define $(y^\mathbb{P},z^\mathbb{P}):=(y^\mathbb{P}(\tau,\xi),z^\mathbb{P}(\tau,\xi))$ as the unique solution of the following standard BSDE (existence and uniqueness have been proved under our assumptions by Kobylanski in \cite{kob})
\begin{equation}\label{bsdep}
y_t^\mathbb P=\xi-\int_t^\tau\widehat F_s(y_s^\mathbb P,z_s^\mathbb P)ds-\int_t^\tau z_s^\mathbb PdB_s, \text{ }0\leq t\leq \tau,\text{ } \mathbb P-a.s.
\end{equation}

First, we introduce the following simple generalization of the comparison theorem proved in \cite{tev} (see Theorem $2$).
\begin{Proposition}\label{prop.comp}
Let Assumptions \ref{assump.hquad} hold true. Let $\xi_1$ and $\xi_2\in L^\infty(\mathbb P)$ for some probability measure $\mathbb P$, and $V^i$, $i=1,2$ be two adapted, c\`adl\`ag non-decreasing processes null at $0$. Let $(Y^i,Z^i)\in \mathbb D^{\infty}(\mathbb P)\times\mathbb H^{2}(\mathbb P)$, $i=1,2$ be the solutions of the following BSDEs
$$Y_t^i=\xi^i-\int_t^T\widehat F_s(Y^i_s,Z^i_s)ds-\int_t^TZ^i_sdB_s+V_T^i-V_t^i,\text{ }\mathbb P-a.s., \text{ }i=1,2,$$
respectively. If $\xi_1\geq \xi_2$, $\mathbb P-a.s.$ and $V^1-V^2$ is non-decreasing, then it holds $\mathbb P-a.s.$ that for all $t\in [0,T]$,
$Y_t^1\geq Y_t^2$.
\end{Proposition}

\proof
First of all, we need to justify the existence of the solutions to those BSDEs. Actually, this is a simple consequence of the existence results of Kobylanski \cite{kob} and for instance Proposition $3.1$ in \cite{ma}. Then, the above comparison is a mere generalization of Theorem $2$ in \cite{tev}.
\ep

\vspace{0.2em}
We then have similarly as in Theorem $4.3$ of \cite{stz} the following results which justifies the PDE intuition given in the Introduction.

\begin{Theorem}\label{quad.unique}
Let Assumptions \ref{assump.hquad} hold. Assume $\xi \in \mathbb{L}^{\infty}_H$ and that $(Y,Z)\in \mathbb D^{\infty}_H\times\mathbb H^{2}_H$ is a solution to $2$BSDE \reff{2bsdequad}. Then, for any $\mathbb{P}\in\mathcal{P}_H$ and $0\leq t_1< t_2\leq T$,
\begin{equation}
Y_{t_1}=\underset{\mathbb{P}^{'}\in\mathcal{P}_H(t_1^+,\mathbb{P})}{\esup^\mathbb{P}}y_{t_1}^{\mathbb{P}^{'}}(t_2,Y_{t_2}), \text{ }\mathbb{P}-a.s.
\label{representationquad}
\end{equation}
Consequently, the $2$BSDE \reff{2bsdequad} has at most one solution in $ \mathbb D^{\infty}_H\times\mathbb H^{2}_H$.
\end{Theorem}

\vspace{0.2em}
Before proceeding with the proof, we will need the following Lemma which shows that in our $2$BSDE framework, we still have a deep link between quadratic growth in $z$ of the generator and the BMO spaces.

\begin{Lemma}\label{lemma.bmo}
Let Assumption \ref{assump.hquad} hold. Assume $\xi \in \mathbb{L}^{\infty}_H$ and that $(Y,Z)\in \mathbb D^{\infty}_H\times\mathbb H^{2}_H$ is a solution to $2$BSDE \reff{2bsdequad}. Then $Z\in\mathbb B\rm{MO}(\mathcal P_H)$.
\end{Lemma}

\vspace{0.2em}
\proof
By It\^o's formula under $\mathbb P$ applied to $e^{-\nu Y_t}$, which is a c\`adl\`ag process, for some $\nu>0$, we have for every $\tau \in \mathcal T^T_0$

\begin{align*}
\frac{\nu^2}{2}\int_\tau^{T}e^{-\nu Y_t}\abs{\widehat a_t^{1/2}Z_t}^2dt&=e^{-\nu \xi}-e^{-\nu Y_{\tau}}-\nu\int_\tau^{T}e^{-\nu Y_{t^-}}dK_t^\mathbb P+\nu\int_\tau^{T}e^{-\nu Y_t}\widehat F_t(Y_t,Z_t)dt\\
&\hspace{0.9em}+\nu\int_\tau^{T}e^{-\nu Y_{t^-}}Z_tdB_t-\sum_{\tau\leq s\leq T}e^{-\nu Y_s}-e^{-\nu Y_{s^{-}}}+\nu\Delta Y_se^{-\nu Y_{s^{-}}}.
\end{align*}

Since $Y\in \mathbb D^{\infty}_H$, $K^\mathbb P$ is non-decreasing and since the contribution of the jumps is negative because of the convexity of the function $x\rightarrow e^{-\nu x}$, we obtain with Assumption \ref{assump.hquad}$\rm{(iv)}$
\begin{align*}
\frac{\nu^2}{2}\mathbb E^\mathbb P_\tau\left[\int_\tau^{T}e^{-\nu Y_t}|\widehat a_t^{\frac12}Z_t|^2dt\right] &\leq e^{\nu\No{Y}_{\mathbb D^{\infty}_H}}{\scriptstyle\left(1+\nu T\left(\alpha+\beta\No{Y}_{\mathbb D^{\infty}_H}\right)\right)}+\frac{\nu\gamma}{2}\mathbb E^\mathbb P_\tau\left[\int_\tau^{T}e^{-\nu Y_t}|\widehat a_t^{\frac12}Z_t|^2dt\right].
\end{align*}

By choosing $\nu=2\gamma$, we then have
$$\mathbb E^\mathbb P_\tau\left[\int_\tau^{T}e^{-2\gamma Y_t}\abs{\widehat a_t^{1/2}Z_t}^2dt\right] \leq \frac{1}{\gamma^2} e^{2\gamma\No{Y}_{\mathbb D^{\infty}_H}}\left(1+2\gamma T\left(\alpha+\beta\No{Y}_{\mathbb D^{\infty}_H}\right)\right).$$

Finally, we obtain
$$\mathbb E^\mathbb P_\tau\left[\int_\tau^{T}\abs{\widehat a_t^{1/2}Z_t}^2dt\right] \leq \frac{1}{\gamma^2} e^{4\gamma\No{Y}_{\mathbb D^{\infty}_H}}\left(1+2\gamma T\left(\alpha+\beta\No{Y}_{\mathbb D^{\infty}_H}\right)\right),$$
which provides the result by arbitrariness of $\mathbb P$ and $\tau$.
\ep

\vspace{0.2em}
{\bf Proof of Theorem \ref{quad.unique}}. The proof follows the lines of the proof of Theorem $4.4$ in \cite{stz}, but we have to deal with some specific difficulties due to our quadratic growth assumption. First if \reff{representationquad} holds, this implies that
$$Y_{t}=\underset{\mathbb{P}^{'}\in\mathcal{P}_H(t^+,\mathbb{P})}{\esup^\mathbb{P}}y_{t}^{\mathbb{P}^{'}}(T,\xi), \text{ } t\in [0,T], \text{ }\mathbb{P}-a.s. \text{ for all }\mathbb P\in \mathcal{P}_H,$$
and thus is unique. Then, since we have that $d\left<Y,B\right>_t=Z_td\left<B\right>_t, \text{ } \mathcal{P}_H-q.s.$, $Z$ is also unique. We now prove \reff{representationquad} in three steps. Roughly speaking, we will obtain one inequality using the comparison theorem, and the second one by using the minimal condition \reff{2bsde.minquad}.

\vspace{0.2em}
$\rm{(i)}$ Fix $0\leq t_1<t_2\leq T$ and $\mathbb P\in\mathcal P_H$. For any $\mathbb P^{'}\in\mathcal P_H(t_1^+,\mathbb P)$, we have
$$Y_{t} = Y_{t_2} - \int_{t}^{t_2} \widehat{F}_s(Y_s,Z_s)ds - \int_{t}^{t_2} Z_sdB_s + K_{t_2}^{ \mathbb P^{'}} - K_{t}^{ \mathbb P^{'} }, \text{ } t_1\leq t\leq t_2, \text{ } \mathbb P^{'}-a.s.$$
and that $K^{\mathbb P^{'}}$ is nondecreasing, $\mathbb P^{'}-a.s.$ Then, we can apply the comparison theorem of Proposition \ref{prop.comp} under $\mathbb P^{'}$ to obtain $Y_{t_1}\geq y_{t_1}^{\mathbb P^{'}}(t_2,Y_{t_2})$, $\mathbb P^{'}-a.s.$ Since $\mathbb P^{'}=\mathbb P$ on $\mathcal F_t^+$, we get $Y_{t_1}\geq y_{t_1}^{\mathbb P^{'}}(t_2,Y_{t_2})$, $\mathbb P-a.s.$ and thus
$$Y_{t_1}\geq\underset{\mathbb{P}^{'}\in\mathcal{P}_H(t_1^+,\mathbb{P})}{\esup^\mathbb{P}}y_{t_1}^{\mathbb{P}^{'}}(t_2,Y_{t_2}), \text{ }\mathbb{P}-a.s.$$

\vspace{0.2em}
$\rm{(ii)}$ We now prove the reverse inequality. Fix $\mathbb P\in\mathcal P_H$. Assume for the time being that
$$C_{t_1}^{\mathbb P,p}:=\underset{\mathbb{P}^{'}\in\mathcal{P}_H(t_1^+,\mathbb{P})}{\esup^\mathbb{P}}\mathbb E_{t_1}^{\mathbb{P}^{'}}\left[\left(K_{t_2}^{\mathbb P^{'}}-K_{t_1}^{\mathbb P^{'}}\right)^p\right]<+\infty,\text{ }\mathbb P-a.s., \text{ for all }p\geq 1.$$

For every $\mathbb P^{'}\in \mathcal P_H(t^+,\mathbb P)$, denote $\delta Y:=Y-y^{\mathbb P^{'}}(t_2,Y_{t_2})\text{ and }\delta Z:=Z-z^{\mathbb P^{'}}(t_2,Y_{t_2}).$ By $\rm{(vi)}$ and $\rm{(v)}$ of Assumption \ref{assump.hquad}, there exist a bounded process $\lambda$ and a process $\eta$ with
$$\abs{\eta_t}\leq \mu\left(\abs{\widehat a_t^{1/2}Z_t}+\abs{\widehat a_t^{1/2}z_t^{\mathbb P^{'}}}\right), \text{ }\mathbb P^{'}-a.s.$$
such that
$$\delta Y_t=\int_t^{t_2}\left(\lambda_s\delta Y_s+(\eta_s+\phi_s)\widehat{a}_s^{1/2}\delta Z_s\right)ds-\int_t^{t_2}\delta Z_sdB_s+K_{t_2}^{\mathbb P^{'}}-K_t^{\mathbb P^{'}}, \text{ }t\leq t_2,\text{ }\mathbb P^{'}-a.s.$$

Define for $t_1\leq t\leq t_2$, $M_t:=\exp\left(\int_{t_1}^t\lambda_sds\right),\text{ }\mathbb P^{'}-a.s.$ Now, since $\phi\in\mathbb B\rm{MO}(\mathcal P_H)$, by Lemma \ref{lemma.bmoholder}, we know that the exponential martingale
$$\mathcal E\left(\int_0^.(\phi_s+\eta_s)\widehat a_s^{-1/2}dB_s\right),$$
is a $\mathbb P^{'}$-uniformly integrable martingale. Therefore we can define a probability measure $\mathbb Q^{'}$, which is equivalent to $\mathbb P^{'}$, by its Radon-Nykodym derivative
$$\frac{d\mathbb Q^{'}}{d\mathbb P^{'}}=\mathcal E\left(\int_0^T(\phi_s+\eta_s))\widehat a_s^{-1/2}dB_s\right).$$

Then, by It\^o's formula, we obtain, as in \cite{stz}, that
$$\delta Y_{t_1}=\mathbb E_{t_1}^{\mathbb Q^{'}}\left[\int_{t_1}^{t_2}M_tdK_t^{\mathbb P^{'}}\right]\leq\mathbb E_{t_1}^{\mathbb Q^{'}}\left[\underset{t_1\leq t\leq t_2}{\sup}(M_t)(K_{t_2}^{\mathbb P^{'}}-K_{t_1}^{\mathbb P^{'}})\right],$$
since $K^{\mathbb P^{'}}$ is non-decreasing. Then, since $\lambda$ is bounded, we have that $M$ is also bounded and thus for every $p\geq 1$
\begin{equation}
\mathbb E_{t_1}^{\mathbb P^{'}}\left[\underset{t_1\leq t\leq t_2}{\sup}(M_t)^p\right]\leq C_p,\text{ } \mathbb P^{'}-a.s.
\label{eq.m}
\end{equation}

This is now that the we will need more work than in the Lipschitz case, and that the BMO properties of the solution will play a crucial role. Since $(\eta + \phi)\widehat a_s^{-1/2}$ is in $\mathbb B\rm{MO}({\mathcal P_H})$, we know by Lemma \ref{lemma.bmoholder} that there exists $r>1$, independent of $\mathbb P$, such that
$$\underset{\mathbb P\in\mathcal P_H}{\Sup}\mathbb E^\mathbb P\left[\left(\mathcal E\left(\int_0^T(\phi_s+\eta_s)\widehat a_s^{-1/2}dB_s\right)\right)^r\right]<+\infty.$$

Then it follows from H\"older inequality (by choosing $r$ as above) and Bayes Theorem that
\begin{align*}
\delta Y_{t_1}&\leq \frac{ \left( \mathbb E^{\mathbb P^{'}}_{t_1} \left[ \mathcal E\left( \int_0^T(\phi_s+\eta_s)\widehat a_s^{-1/2}dB_s\right)^r\right ]\right )^{\frac1r}}{\mathbb E^{\mathbb P^{'}}_{t_1}\left[\mathcal E\left(\int_0^T(\phi_s+\eta_s)\widehat a_s^{-1/2}dB_s\right)\right]}\left(\mathbb E_{t_1}^{\mathbb P^{'}}\left[\left(\underset{t_1\leq t\leq t_2}{\sup}M_t\right)^q\left(K_{t_2}^{\mathbb P^{'}}-K_{t_1}^{\mathbb P^{'}}\right)^q\right]\right)^{\frac1q}\\
&\leq C\left(\mathbb E_{t_1}^{\mathbb P^{'}}\left[\left(K_{t_2}^{\mathbb P^{'}}-K_{t_1}^{\mathbb P^{'}}\right)^{2q-\frac12}\left(K_{t_2}^{\mathbb P^{'}}-K_{t_1}^{\mathbb P^{'}}\right)^{\frac12}\right]\right)^{\frac{1}{2q}}\\
&\leq C \left( C_{t_1}^{\mathbb P, 4q-1}\right)^{\frac{1}{4q}} \left( \mathbb E_{t_1}^{\mathbb P^{'}}\left[K_{t_2}^{\mathbb P^{'}}-K_{t_1}^{\mathbb P^{'}}\right]\right)^{\frac{1}{4q}}.
\end{align*}

By the minimum condition \reff{2bsde.minquad} and since $\mathbb P^{'}\in \mathcal P_H(t^+,\mathbb P)$ is arbitrary, this ends the proof.

\vspace{0.2em}
$\rm{(iii)}$ It remains to show that the estimate for $C_{t_1}^{\mathbb P,p}$ holds for $p\geq 1$. Once again, this will be possible because of the BMO property satisfied by the solution of the 2BSDE. By definition of the family $\left\{K^\mathbb P, \mathbb P\in\mathcal P_H\right\}$, we have
\begin{align*}
\mathbb E^{\mathbb{P}^{'}}\left[\left(K_{t_2}^{\mathbb P^{'}}-K_{t_1}^{\mathbb P^{'}}\right)^p\right] &\leq C\left({\scriptstyle1+\No{Y}^p_{\mathbb D^{\infty}_H}+\No{\xi}^p_{\mathbb L^{\infty}_H}}+\mathbb E_{t_1}^{\mathbb{P}^{'}}\left[\left(\int_{t_1}^{t_2}|\widehat a^{\frac12}_tZ_t|^2dt\right)^p+\left(\int_{t_1}^{t_2}Z_tdB_t\right)^p\right]\right).
\end{align*}

Thus by the energy inequalities \reff{energy} and by Burkholder-Davis-Gundy inequality, we obtain
$$\mathbb E^{\mathbb{P}^{'}}\left[\left(K_{t_2}^{\mathbb P^{'}}-K_{t_1}^{\mathbb P^{'}}\right)^p\right]\leq C\left(1+\No{Y}^p_{\mathbb D^{\infty}_H}+\No{\xi}^p_{\mathbb L^{\infty}_H}+\No{Z}^{2p}_{\mathbb B\rm{MO}_H}+\No{Z}^{p}_{\mathbb B\rm{MO}_H}\right)<+\infty.$$

Then we can proceed exactly as in the proof of Theorem $4.3$ in \cite{stz}.
\ep

\vspace{0.2em}
\begin{Remark}
It is interesting to notice that in contrast with standard quadratic BSDEs, for which the only property of BMO martingales used to obtain uniqueness is the fact that their Dol\'eans-Dade exponential is a uniformly integrable martingale, we need a lot more in the 2BSDE framework. Indeed, we use extensively the energy inequalities and the existence of moments for the Dol\'eans-Dade exponential (which is a consequence of the so called reverse H\"older inequalities, which is a more general version of Lemma \ref{lemma.bmoholder}). Furthermore, we will also use the so-called Muckenhoupt condition (which corresponds to Lemma \ref{lemma.bmo2}, see \cite{kaz} for more details) in both our proofs of existence. This seems to be directly linked to the presence of the non-decreasing processes $K^\mathbb P$ and raises the question about whether it could be possible to generalize the recent approach of Barrieu and El Karoui \cite{elkarbar}, to second-order BSDEs. Indeed, since they no longer assume a bounded terminal condition, the $Z$ part of the solution is no-longer BMO. We leave this interesting but difficult question to future research.
\end{Remark}

We conclude this section by showing some a priori estimates which will be useful in the sequel. Notice that these estimates also imply uniqueness, but they use intensively the representation formula \reff{representationquad}.

\begin{Theorem}\label{estimatesquad}
Let Assumption \ref{assump.hquad} hold.

\vspace{0.2em}
(i) Assume that $\xi\in\mathbb L^{\infty}_H$ and that $(Y,Z)\in \mathbb D^{\infty}_H\times\mathbb H^{2}_H$ is a solution to $2$BSDE \reff{2bsdequad}. Then, there exists a constant $C$ such that for all $p\geq 1$
\begin{align*}
&\No{Y}_{\mathbb D^{\infty}_H}+\No{Z}_{\mathbb B\rm{MO}(\mathcal P_H)}^2\leq C\left(1+\No{\xi}_{\mathbb L^{\infty}_H}\right),\ \underset{\mathbb P\in \mathcal P_H,\text{ }\tau\in\mathcal T_0^T}{\sup}\mathbb E_\tau^\mathbb P\left[(K_T^\mathbb P-K_\tau^\mathbb P)^p\right]\leq C \left(1+\No{\xi}_{\mathbb L^{\infty}_H}^p\right).
\end{align*}

\vspace{0.2em}
(ii) Assume that $\xi^i\in\mathbb L^{\infty}_H$ and that $(Y^i,Z^i)\in \mathbb D^{\infty}_H\times\mathbb H^{2}_H$ is a corresponding solution to $2$BSDE \reff{2bsdequad}, $i=1,2$. Denote $\delta\xi:=\xi^1-\xi^2$, $\delta Y:=Y^1-Y^2$, $\delta Z:=Z^1-Z^2$ and $\delta K^\mathbb P:=K^{\mathbb P,1}-K^{\mathbb P,2}$. Then, there exists a constant $C$ such that
\begin{align*}
&\No{\delta Y}_{\mathbb D^{\infty}_H}\leq C\No{\delta\xi}_{\mathbb L^{\infty}_H},\
\No{\delta Z}_{\mathbb B\rm{MO}(\mathcal P_H)}^2\leq C\No{\delta\xi}_{\mathbb L^{\infty}_H}\left(1+\No{\xi^1}_{\mathbb L^{\infty}_H}+\No{\xi^2}_{\mathbb L^{\infty}_H}\right)\\
& \hspace{0.5cm} \forall p\geq 1, \text{ } \underset{\mathbb P\in \mathcal P_H}{\sup}\mathbb E^\mathbb P\left[\underset{0\leq t\leq T}{\sup}\abs{\delta K_t^\mathbb P}^p\right]\leq C \No{\xi}_{\mathbb L^{\infty}_H}^{p/2}\left(1+\No{\xi^1}_{\mathbb L^{\infty}_H}^{p/2}+\No{\xi^2}_{\mathbb L^{\infty}_H}^{p/2}\right).
\end{align*}
\end{Theorem}

\proof
$\rm{(i)}$ By Lemma $1$ in \cite{bh}, for all $\mathbb P\in \mathcal P_H$, $|y_t^{\mathbb P}|\leq \alpha\left(e^{\beta T}-1\right)/\beta+ e^{\beta T}\No{\xi}_{\mathbb L^{\infty}_H},$ and by \reff{representationquad}, the estimate of $\No{Y}_{\mathbb D^{\infty}_H}$ is clear. By the proof of Lemma \ref{lemma.bmo}, we have
\begin{align*}
\No{Z}_{\mathbb B\rm{MO}(\mathcal P_H)}^2 &\leq Ce^{C\No{Y}_{\mathbb D^{\infty}_H}}\left(1+\No{Y}_{\mathbb D^{\infty}_H}\right)\leq C\left(1+\No{\xi}_{\mathbb L^{\infty}_H}\right).
\end{align*}

Finally, we have for all $\tau\in\mathcal T^T_0$, for all $\mathbb P\in \mathcal P_H$ and for all $p\geq 1$, by definition
$$(K_T^\mathbb P-K_\tau^\mathbb P)^p=\left(Y_\tau-\xi+\int_\tau^T\widehat F_t(Y_y,Z_t)dt+\int_\tau^TZ_tdB_t\right)^p.$$

Therefore, by our growth Assumption \ref{assump.hquad}$\rm{(iv)}$
\begin{align*}
\mathbb E^\mathbb P_\tau\left[(K_T^\mathbb P-K_\tau^\mathbb P)^p\right] & \leq C\left(1+\No{\xi}_{\mathbb L^{\infty}_H}^p+\No{Y}_{\mathbb D^{\infty}_H}^p+\mathbb E^\mathbb P_\tau\left[\left(\int_{\tau}^{T}|\widehat a^{\frac12}_tZ_t|^2dt\right)^p+\left(\int_{\tau}^{T}Z_tdB_t\right)^p\right]\right)\\
&\leq C\left(1+\No{\xi}_{\mathbb L^{\infty}_H}^p+\No{Z}_{\mathbb B\rm{MO}(\mathcal P_H)}^{2p}+\No{Z}_{\mathbb B\rm{MO}(\mathcal P_H)}^{p}\right)\leq C\left(1+\No{\xi}_{\mathbb L^{\infty}_H}^p\right),
\end{align*}
where we used again the energy inequalities and the BDG inequality. This provides the estimate for $K^\mathbb P$ by arbitrariness of $\tau$ and $\mathbb P$.

\vspace{0.2em}
$\rm{(ii)}$ With the same notations and calculations as in step $\rm{(ii)}$ of the proof of Theorem \ref{quad.unique}, it is easy to see that for all $\mathbb P\in \mathcal P_H$ and for all $t\in[0,T]$, we have $\delta y_t^\mathbb P=\mathbb E^\mathbb Q_t\left[M_T\delta\xi\right]\leq C\No{\delta\xi}_{\mathbb L^{\infty}_H},$ since $M$ is bounded and we have \reff{eq.m}. By Theorem \ref{quad.unique}, the estimate for $\delta Y$ follows immediately, since we have for all $\mathbb P\in\mathcal P_H$
$$\abs{\delta Y_t}\leq \underset{\mathbb P^{'}\in\mathcal P_H(t^+,\mathbb P)}{\esup^\mathbb P}\abs{\delta y^{\mathbb P^{'}}},\ \mathbb P-a.s.$$

Denote $\delta \widehat F^{1,2}_t:=\widehat F_t(Y_t^1,Z_t^1)-\widehat F_t(Y_t^2,Z_t^2)$. Apply It\^o's formula to $\abs{\delta Y}^2$ for $\tau\in\mathcal T_0^T$ and $T$
\begin{align*}
\mathbb E^\mathbb P_\tau\left[\abs{\delta Y_\tau}^2+\int_{\tau}^{T}|\widehat a^{\frac12}_t\delta Z_t|^2dt\right]&\leq\mathbb E^\mathbb P_\tau\left[\abs{\delta \xi}^2-2\int_{\tau}^{T}\delta Y_t\delta\widehat F^{1,2}_tdt-2\int_\tau^T\delta Y_{t^-}d(\delta K_t^\mathbb P)\right].
\end{align*}

Then, we have by Assumption \ref{assump.hquad}$\rm{(iv)}$ and the estimates proved in $\rm{(i)}$ above
\begin{align*}
\mathbb E^\mathbb P_\tau\left[\int_{\tau}^{T}\abs{\widehat a^{1/2}_t\delta Z_t}^2dt\right]&\leq C\No{\delta Y}_{\mathbb D^{\infty}_H}\left(1+\sum_{i=1}^2\No{Y^i}_{\mathbb D^{\infty}_H}+ \No{Z^i}_{\mathbb B\rm{MO}(\mathcal P_H)}\right)\\
&+\No{\delta\xi}_{\mathbb L^{\infty}_H}^2+2\No{\delta Y}_{\mathbb D^{\infty}_H}\mathbb E^\mathbb P_\tau\left[\abs{K_T^{\mathbb P,1}-K_\tau^{\mathbb P,1}}+\abs{K_T^{\mathbb P,2}-K_\tau^{\mathbb P,2}}\right]\\
&\leq C\No{\delta\xi}_{\mathbb L^{\infty}_H}\left(1+\No{\xi^1}_{\mathbb L^{\infty}_H}+\No{\xi^2}_{\mathbb L^{\infty}_H}\right),
\end{align*}
which implies the required estimate for $\delta Z$. Finally, by definition, we have for all $\mathbb P\in\mathcal P_H$
$$\delta K^\mathbb P_t=\delta Y_0-\delta Y_t-\int_0^t\left[\widehat F_s(Y^1_s,Z^1_s)-\widehat F_s(Y^2_s,Z^2_s)\right]ds +\int_0^t\delta Z_sdB_s,\ t\in[0,T],\ \mathbb P-a.s.$$

By Assumptions \ref{assump.hquad}$\rm{(v)}$ and $\rm{(vi)}$, it follows that
\begin{align*}
\underset{0\leq t\leq T}{\sup} \abs{\delta K_t^\mathbb P} &\leq C\left( \No{\delta Y}_{\mathbb D^{\infty}_H}+\int_0^T|\widehat a_s^{\frac12}\delta Z_s|(1+|\widehat a_s^{\frac12} Z_s^1|+|\widehat a_s^{\frac12}Z_s^2|)ds+\underset{0\leq t\leq T}{\sup}\abs{\int_0^t\delta Z_sdB_s}\right),
\end{align*}

and by Cauchy-Schwarz, BDG and energy inequalities, we see that
\begin{align*}
\mathbb E^\mathbb P\left[\underset{0\leq t\leq T}{\sup} \abs{\delta K_t^\mathbb P}^p\right] &\leq C\mathbb E^\mathbb P\left[\left(\int_0^T1+|\widehat a_s^{\frac12} Z_s^1|^2+|\widehat a_s^{\frac12}Z_s^2|^2ds\right)^p\right]^{\frac12} \mathbb E^\mathbb P\left[\left(\int_0^T|\widehat a_s^{\frac12}\delta Z_s|^2ds\right)^p\right]^{\frac12}\\
&\hspace{0.9em}+C\left( \No{\delta \xi}_{\mathbb L^{\infty}_H}^p+\mathbb E^\mathbb P\left[\left(\int_0^T\abs{\widehat a_s^{1/2}\delta Z_s}^2ds\right)^{p/2}\right]\right)\\
&\leq C\No{\delta \xi}_{\mathbb L^{\infty}_H}^{p/2}\left(1+\No{ \xi^1}_{\mathbb L^{\infty}_H}^{p/2}+\No{ \xi^2}_{\mathbb L^{\infty}_H}^{p/2}\right).
\end{align*}
\ep

\begin{Remark}\label{remark.proofquad}
Let us note that the proof of $\rm{(i)}$ only requires that Assumption \ref{assump.hquad}$\rm{(iv)}$ holds true, whereas $\rm{(ii)}$ also requires Assumption \ref{assump.hquad}$\rm{(v)}$ and $\rm{(vi)}$.
\end{Remark}

\section{Purely quadratic 2BSDEs and robust risk-sensitive control}\label{sec.pure}
This section is devoted to simple examples of quadratic 2BSDEs for which existence can be obtained in a straightforward way. It is our conviction that starting from these easy situations will allow the reader to realize that, fortunately, some techniques from the classical quadratic BSDE literature still work in our framework, albeit with more complications (conversely, the next Section will explain why from our point of view, 
it's a very difficult task to try and obtain general versions of the monotone approximation techniques). In this regard, we will first introduce the "simplest" quadratic 2BSDEs corresponding to a so-called purely quadratic generator. Keeping in mind the possible applications of this theory, we will then define a quasi-sure version of the entropic risk measure which we will then use to study robust version study risk-sensitive control problems.

\subsection{Entropy and purely quadratic 2BSDEs}\label{purely}
Given $\xi\in \mathcal L^\infty_H$, we first consider the purely quadratic 2BSDE defined as follows
\begin{equation}
Y_t=-\xi +\int_t^T\frac{\gamma}{2}\abs{\widehat{a}_s^{1/2}Z_s}^2ds -\int_t^TZ_sdB_s + K_T^\mathbb P-K_t^\mathbb P, \text{ } 0\leq t\leq T, \text{  } \mathcal P_H-q.s.
\label{2bsdepurequad}
\end{equation}
Then we use the classical exponential change of variables and define 
$$\overline{Y}_t:=e^{\gamma Y_t}, \text{ }\overline{Z}_t:=\gamma \overline{Y}_tZ_t, \text{ }\overline{K}^{\mathbb P}_t:=\gamma\int_0^t\overline{Y}_sdK^{\mathbb P}_s-\sum_{0\leq s\leq t}e^{\gamma Y_s}-e^{\gamma Y_{s^-}}-\gamma\Delta Y_s e^{\gamma Y_{s^-}}.$$
At least formally, we see that $(\overline{Y},\overline{Z},\overline{K}^{\mathbb P})$ verifies the following equation
\begin{equation}
\overline{Y}_t=e^{-\gamma\xi} -\int_t^T\overline{Z}_sdB_s + \overline{K}^{\mathbb P}_T-\overline{K}^{\mathbb P}_t, \text{ } 0\leq t\leq T, \text{  } \mathbb P-a.s.\ \forall\mathbb P\in\mathcal P_H
\label{2bsdepurelip}
\end{equation}
which is in fact a 2BSDE with generator equal to $0$ (and thus Lipschitz), provided that the family $\left(\overline{K}^{\mathbb P}\right)_{\mathbb P\in\mathcal P_H}$ satisfies the minimum condition \reff{2bsde.minquad}. Thus the purely quadratic 2BSDE \reff{2bsdepurequad} is linked to the 2BSDE with Lipschitz generator (\ref{2bsdepurelip}), which has a unique solution by Soner, Touzi and Zhang \cite{stz}. We now make this rigorous.

\begin{Proposition}\label{prop.pure}
The 2BSDE \reff{2bsdepurequad} has a unique solution $(Y,Z)\in\mathbb D^\infty_H\times\mathbb H^2_H$ given by
$$Y_t=\frac1\gamma\ln\left(\underset{\mathbb P^{'}\in\mathcal P_H(t^+,\mathbb P)}{\esup^\mathbb P}\mathbb E^{\mathbb P^{'}}_t\left[e^{-\gamma \xi}\right]\right), \ \mathbb P-a.s., \ t\in[0,T], \text{ for all $\mathbb P\in\mathcal P_H$}.$$
\end{Proposition}

We emphasize that unlike in the classical case where the above result is an easy consequence of It\^o's formula, in the second order framework the situation is complicated by the presence of the non-decreasing processes, and the main difficulty is to show that they do satisfy the required minimum condition.

\vspace{0.2em}
\proof
Uniqueness is a simple consequence of Theorem \ref{quad.unique}. In the following, we prove the existence in 3 steps.

\vspace{0.2em}
{\bf{Step $1$}}: Let $(\overline{Y},\overline{Z})\in\mathbb D^2_H\times\mathbb H^2_H$ be the unique solution to the 2BSDE (\ref{2bsdepurelip}) and $\overline{K}^{\mathbb P}$ be the corresponding non-decreasing processes. In particular, we know that $$\overline{Y}_t=\underset{\mathbb P^{'}\in\mathcal{P}_{H}(t^+,\mathbb P)}{\esup^{\mathbb P}}\mathbb E^{\mathbb P^{'}}_t\left[e^{-\gamma\xi}\right],\ \mathbb P-a.s.,$$
which implies that $\overline{Y}\in\mathbb D^\infty_H$, since $0<e^{-\gamma\No{\xi}_{\mathbb L^\infty_H}}\leq Y_t\leq e^{\gamma\No{\xi}_{\mathbb L^\infty_H}}.$ We can therefore make the following change of variables
$$Y_{t}:=\frac{1}{\gamma}\text{ln}\left( \overline{Y}_t\right),\ Z_{t}:=\frac{1}{\gamma}\frac{\overline{Z}_t}{\overline{Y}_t}, \ {K}^{\mathbb P}_t:=\int^t_0\frac{1}{\gamma \overline{Y}_s}d\overline{K}^{\mathbb P,c}_s-\sum_{0<s\leq t}\frac{1}{\gamma}\text{log}\left( 1-\frac{\Delta \overline{K}^{\mathbb P,d}_s}{\overline{Y}_{s-}} \right).$$
Then by It\^o's formula, we can verify that the triplet $(Y,Z,K^\mathbb P)$ satisfies (\ref{2bsdepurequad}). Moreover, notice that ${K}^{\mathbb P}$ is non-decreasing with ${K}^{\mathbb P}_0=0$.

\vspace{0.2em}
{\bf{Step $2$}}: Denote now $(y^{\mathbb P},z^{\mathbb P})$ the solutions of the standard BSDEs corresponding to the 2BSDE \reff{2bsdepurequad} (existence and uniqueness are ensured for example by \cite{kob}). Furthermore, if we define $\overline{y}^{\mathbb P}_t:=e^{\gamma y^{\mathbb P}_t}, \text{ }\overline{z}^{\mathbb P}_t:=\gamma \overline{y}^{\mathbb P}_tz^{\mathbb P}_t,$
then we know that $(\overline{y}^{\mathbb P},\overline{z}^{\mathbb P})$ solve the standard BSDE under $\mathbb P$ corresponding to \reff{2bsdepurelip}. Due to the monotonicity of the function $x\rightarrow\ln(x)$ and the representation for $\overline{Y}$, we have the following representation for $Y$ $$Y_t=\underset{\mathbb P^{'}\in\mathcal{P}_{H}(t^+,\mathbb P)}{\esup^{\mathbb P}}y^{\mathbb P}_t=\frac1\gamma\ln\left(\underset{\mathbb P^{'}\in\mathcal P_H(t^+,\mathbb P)}{\esup^\mathbb P}\mathbb E^{\mathbb P^{'}}_t\left[e^{-\gamma \xi}\right]\right), \ \mathbb P-a.s.$$

\vspace{0.2em}
{\bf{Step $3$}}: Finally, it remains to check the minimum condition for the family of non-decreasing processes $\left\{K^{\mathbb P}\right\}$. Let $\mathbb P\in\mathcal P_H$, $t\in[0,T]$ and $\mathbb P^{Ô}\in\mathcal P_H(t^+,\mathbb P)$. Using the notations of the proof of Theorem \ref{quad.unique}, denote $\mathcal E_t:=\mathcal E\left(\int_0^t(\phi_s+\eta_s)\widehat a_s^{-1/2}dB_s\right)$. Since the purely quadratic generator satisfies the Assumption \ref{assump.hquad}, we know from this same proof
\begin{align*}
\nonumber \delta Y_{t}=\mathbb E_{t}^{\mathbb Q^{'}}\left[\int_{t}^{T}M_tdK_t^{\mathbb P^{'}}\right]&\geq\mathbb E_{t}^{\mathbb Q^{'}}\left[\underset{t\leq s\leq T}{\inf}(M_s)(K_{T}^{\mathbb P^{'}}-K_{t}^{\mathbb P^{'}})\right]= \mathbb E_{t}^{\mathbb P^{'}}\left[\frac{\mathcal E_T}{\mathcal E_t}\underset{t\leq s\leq T}{\inf}(M_s)(K_{T}^{\mathbb P^{'}}-K_{t}^{\mathbb P^{'}})\right].
\label{remarque}
\end{align*}

Denote $\delta K^{\mathbb P^{'}}:=K_T^{\mathbb P^{'}}-K_t^{\mathbb P^{'}}$. Let $r$ be given by Lemma \ref{lemma.bmo2} applied to $\mathcal E$. Then we estimate
\begin{align*}
&\mathbb E_t^{\mathbb P^{'}}\left[\delta K^{\mathbb P^{'}}\right]\leq \mathbb E_t^{\mathbb P^{'}}\left[\frac{\mathcal E_T}{\mathcal E_t}\underset{t\leq s\leq T}{\inf}(M_s)\delta K^{\mathbb P^{'}}\right]^{\frac{1}{2r-1}} \mathbb E_t^{\mathbb P^{'}}\left[\left(\frac{\mathcal E_T}{\mathcal E_t}\underset{t\leq s\leq T}{\inf}(M_s)^{-1}\right)^{\frac{1}{2(r-1)}}\delta K^{\mathbb P^{'}}\right]^{\frac{2(r-1)}{2r-1}}\\
&\leq\left(\delta Y_t\right)^{\frac{1}{2r-1}}\left(\mathbb E_t^{\mathbb P^{'}}\left[\left(\frac{\mathcal E_T}{\mathcal E_t}\right)^{\frac{1}{r-1}}\right]\right)^{\frac{r-1}{2r-1}}\left(\mathbb E_t^{\mathbb P^{'}}\left[\underset{t\leq s\leq T}{\inf}(M_s)^{-\frac{2}{r-1}}\right]\mathbb E_t^{\mathbb P^{'}}\left[\left(\delta K^{\mathbb P^{'}}\right)^4\right]\right)^{\frac{r-1}{2(2r-1)}}\\
&\leq C\left(\mathbb E_t^{\mathbb P^{'}}\left[\left(K_{T}^{\mathbb P^{'}}\right)^4\right]\right)^{\frac{r-1}{2(2r-1)}}\left(\delta Y_t\right)^{\frac{1}{2r-1}}.
\end{align*}

By following the arguments of the proof of Theorem \ref{quad.unique} $\rm{(ii)}$ and $\rm{(iii)}$, we then deduce the minimum condition.
\ep

\begin{Remark}\label{remark.quad}
In order to prove the minimum condition it is fundamental that the process $M$ above is bounded from below. For instance, it would not be the case if we had replaced the Lipschitz assumption on $y$ by a monotonicity condition as in \cite{pos}.
\end{Remark}

\vspace{0.2em}
\begin{Remark}\label{rem.exist}
More generally, by the same exponential change and arguments above, we can also prove that there exists a unique solution to 2BSDEs with terminal condition $\xi\in \mathcal L^\infty_H$ and generators of the type $\widehat{a}^{1/2}zg(t,\omega)+h(t,\omega)+\frac{\theta}{2}\abs{\widehat{a}^{1/2}_tz}^2 $ where $g$ and $h$ are assumed to be bounded, adapted and uniformly continuous in $\omega$ for the $\No{\cdot}_{\infty} $. This will be useful in the next subsection.
\end{Remark}

Thanks to the above result, we can define a quasi-sure (or robust) version of the dynamic entropic risk measure under volatility uncertainty
$$e_{\gamma,t}(\xi):=\frac1\gamma\ln\left(\underset{\mathbb P^{'}\in\mathcal P_H(t^+,\mathbb P)}{\esup^\mathbb P}\mathbb E^{\mathbb P^{'}}_t\left[e^{-\gamma \xi}\right]\right)=\frac1\gamma\ln\left(\underset{\mathbb P^{'}\in\mathcal P_H(t,\mathbb P)}{\esup^\mathbb P}\mathbb E^{\mathbb P^{'}}_t\left[e^{-\gamma \xi}\right]\right),$$
where the parameter $\gamma$ stands for the risk tolerance and where we used the fact that, as proved in \cite{stz2} (see Proposition $4.11$), the solution of \reff{2bsdepurequad} is actually $\mathbb F$-measurable. This result is a non trivial (albeit expected) extension of the results obtained by El Karoui and Barrieu \cite{bar} on the representation of the dynamic entropic risk measure as the unique solution of a purely quadratic BSDE.


\subsection{An application to robust risk-sensitive control}\label{sec.ex}
One of the possible applications of quadratic BSDEs is to study risk-sensitive control problems. In a classical setting we refer the reader to \cite{bie}, \cite{flem00}, \cite{flem02} or \cite{elkarmat} and the references therein for more details. In this section, we will define and solve a robust version of these problems. First of all, for technical reasons, we restrict the probability measures in $\widetilde{\mathcal{P}}_H:=\widetilde{\mathcal{P}}_S\bigcap\mathcal P_H$, where $\widetilde{\mathcal{P}}_S$ is defined in Subsection \ref{subsec.proba}. Then $\widehat{a}$ is uniformly bounded. For each $\mathbb P\in\widetilde{\mathcal P}_H$, we can define a $\mathbb P$-Brownian motion $W^\mathbb P$ by $dW^{\mathbb P}_t=\widehat{a}_t^{-1/2}dB_t\ \mathbb P-a.s. $ Let us now consider some system, whose evolution is described (for simplicity) by the canonical process $B$. A controller then intervenes on the system via an adapted stochastic process $u$ which takes its values in a compact metric space $U$. The set of those controls is called admissible and denoted by $\mathcal U$. When the controller acts with $u$ under the probability $\mathbb P\in\widetilde{\mathcal P}_H$, the dynamic of the controlled system remains the same, but now under the probability measure $\mathbb P^u$ defined by its density with respect to $\mathbb P$
$$\frac{d\mathbb P^u}{d\mathbb P}=\text{exp}\left(\int_0^T\widehat{a}_t^{-1/2}g(t,B.,u_t)dW^{\mathbb P}_t-\frac{1}{2}\int_0^T\left|\widehat{a}_t^{-1/2}g(t,B.,u_t)\right|^2dt\right), $$
where $g(t,\omega,u)$ is assumed to be bounded uniformly in $u$, continuous with respect to $u$, adapted and uniformly continuous in $\omega$. Notice that this probability measure is well defined
since $\widehat{a}$ is uniformly bounded. Then, under $\mathbb P^u$, the dynamic of the system is given by
$$dB_t=g(t,B.,u_t)dt+\widehat{a}_t^{1/2}dW^{\mathbb P,u},\ \mathbb P^u-a.s. $$
where $W^{\mathbb P,u}$ is a Brownian motion under $\mathbb P^u$ defined by $dW_t^{\mathbb P,u}=dW_t^\mathbb P-\widehat a^{-1/2}_tg(t,B_.,u_t)dt.$ When the controller is risk averse, we assume that the reward functional of the control action is given by the following expression
$$\forall u\in\mathcal U,\,J(u):=\underset{\mathbb P\in\widetilde{\mathcal P}_H}{\sup} \mathbb E^{\mathbb P^u}\left[\exp\left(\theta\int^T_0h(s,B.,u_s)ds+\Psi(B_T)\right)\right],$$
where $\theta>0$ is a real parameter which represents the sensitiveness of the controller with respect to risk. Here $h(t,\omega,u)$ is assumed to be adapted and continuous in $u$, and both $\Psi$ and $h$ are assumed to be bounded and uniformly continuous in $\omega$ for the $\No{\cdot}_{\infty}$ norm. We are interested in finding an admissible control $u^{\ast}$ which maximizes the reward $J(u)$ for the controller. We begin with establishing the link between $J(u)$ and 2BSDEs in the following proposition
\begin{Proposition}\label{quad.u}
There exists a unique solution $(Y^u,Z^u)$ of the 2BSDE associated with the generator
$\widehat a_t^{1/2}zg(t,B.,u_t)+h(t,B.,u_t)+\frac{\theta}{2}|\widehat{a}^{1/2}_tz|^2 $, \textit{i.e.}, $ \mathbb P-a.s.,\text{ for all $\mathbb P\in\widetilde{\mathcal P}_H$}$
\begin{equation}\label{sen.ex}
Y^u_t=\Psi(B_T)+\int_t^T\left(\widehat a_s^{1/2}Z^u_sg(s,B.,u_s)+h(s,B.,u_s)+\frac{\theta}{2}|\widehat{a}^{1/2}_sZ^u_s|^2\right)ds-\int_t^TZ^u_sdB_s-dK^{u,\mathbb P}_t.
\end{equation}
Moreover $J(u)=\exp\left(\theta Y^u_0\right)$.
\end{Proposition}

\proof
With our assumptions on $g$, $h$ and $\Psi$, we know from Remark \ref{rem.exist} that there exists a unique solution to the 2BSDE (\ref{sen.ex}). Using the results of \cite{elkarmat} in the non-robust case and by the representation for $Y^u$, we have
$$Y^{u}_t=\frac{1}{\theta}\ln\left(\underset{\mathbb{P'}\in\widetilde{\mathcal P}_H(t^+,\mathbb P)}{\esup^{\mathbb P}}\mathbb E^{\mathbb{P'}^u}_t\left[\exp\left(\theta\int^T_th(s,B.,u_s)ds+\Psi(B_T)\right)\right]\right),\ \mathbb P-a.s.$$
Therefore, we have $J(u)=\text{exp}\left\{\theta Y^u_0\right\}$.
\ep

\vspace{0.2em}
As explained in \cite{elkarmat}, by applying Benes' selection theorem, there exists a measurable version $u^{\ast}(t,B.,z)$ of $\text{arg max}\ I(t,B.,z,u):=\widehat a_t^{1/2}(B_.)zg(t,B.,u)+h(t,B.,u).$ We know that $I^{\ast}(t,B.,z):={\sup}_{u\in U}I(t,B.,z,u)=I(t,B.,z,u^{\ast}(t,B.,z))$ is convex uniformly Lipschitz in $z$ because it is the supremum of functions which are linear in $z$. So the mapping $z \rightarrow I^{\ast}(t,B.,z)+\frac{1}{2}|\widehat{a}^{1/2}_tz|^2$ is continuous with quadratic growth, implying that a solution $(y^{\ast,\mathbb P},z^{\ast,\mathbb P})$ of the BSDE associated to this generator exists.
Then we have
\begin{Theorem}
There exists a unique solution $(Y^{\ast},Z^{\ast}) $ to the following 2BSDE
\begin{equation}\label{quad.ast}
Y^{\ast}_t=\Psi(B_T)+\int_t^T\left(I^{\ast}(s,B.,Z^{\ast}_s)+\frac{\theta}{2}|\widehat{a}^{1/2}_sZ^{\ast}_s|^2\right)ds-\int_t^TZ^{\ast}_sdB_s+K^{\ast,\mathbb P}_T-K^{\ast,\mathbb P}_t.
\end{equation}

The admissible control $u^{\ast}:=(u^{\ast}(t,B.,Z^{\ast}_t))_{t\leq T} $ is optimal and $(\exp(\theta Y^{\ast}_t))_{t\leq T} $ is the value function of the robust risk-sensitive control problem, i.e., for any $t\leq T$ we have:
\begin{equation}\label{optimale}
\exp(\theta Y^{\ast}_t) =\underset{\mathbb{P'}\in\widetilde{\mathcal P}_H(t^+,\mathbb P)}{\esup^{\mathbb P}} \underset{u\in\mathcal U}{\esup^{\mathbb P}}\ \mathbb E^{\mathbb {P'}^{u}}_t\left[\exp\left(\theta\int^T_th(s,B.,u_s)ds+\Psi(B_T)\right)\right].
\end{equation}
\end{Theorem}

\vspace{0.2em}
\proof
First, we need to prove the existence of a solution to the quadratic 2BSDE (\ref{quad.ast}). Unlike in Proposition \ref{quad.u}, here $u^{\ast}$ also depends on $z$, so we are not directly in the context of Remark \ref{rem.exist}. However, it is easy to see that the generator satisfies Assumption \ref{assump.hquad}, which ensures that we have uniqueness of the solution. Therefore, exactly as in Proposition \ref{quad.u}, for $\mathbb P\in\widetilde{\mathcal P}_H$, by making the exponential change
$$\overline{Y}_t:=e^{\theta Y^{\ast}_t}, \text{ }\overline{Z}_t:=\theta \overline{Y}_tZ^{\ast}_t, \text{ }\overline{K}^{\mathbb P}_t:=\theta\int_0^t\overline{Y}_sdK^{\ast,\mathbb P}_s-\sum_{0\leq s\leq t}e^{\theta Y^{\ast}_s}-e^{\theta Y^{\ast}_{s^-}}-\theta\Delta Y^{\ast}_s e^{\theta Y^{\ast}_{s^-}},$$
we see that $(\overline{Y},\overline{Z},\overline{K}^{\mathbb P})$ formally verifies the following equation, $\mathbb P-a.s.$ for every $\mathbb P\in\widetilde{\mathcal P}_H$
\begin{equation}
\overline{Y}_t=e^{\theta\Psi(B_T)} +\int_t^T \underset{u\in U}{\sup}\left\{\widehat a_s^{1/2}\overline{Z}_sg(s,B.,u)+\theta \overline{Y}_sh(s,B.,u)\right\}ds -\int_t^T\overline{Z}_sdB_s + \overline{K}^{\mathbb P}_T-\overline{K}^{\mathbb P}_t.
\label{2bsdequalip}
\end{equation}

Since this is 2BSDE with Lipschitz generator from Soner, Touzi and Zhang \cite{stz}, we know that $(\overline{Y},\overline{Z},\overline{K}^{\mathbb P})$ exists, is unique and satisfies the representation property \reff{representationquad}. Arguing exactly as in Subsection \ref{purely} for the purely quadratic 2BSDEs, we can then obtain the existence. Now, from \cite{elkarmat}, we have that $$\exp\left(\theta y^{\ast,\mathbb P}_t\right) =\underset{u\in\mathcal U}{\esup^{\mathbb P}}\ \mathbb E^{\mathbb P^u}_t\left[\exp\left(\theta\int^T_th(s,B.,u_s)ds+\Psi(B_T)\right)\right].$$
Then the representation for $Y^{\ast} $ implies the desired result.
\ep

\begin{Remark}
We acknowledge that a robust problem for the controller may be more naturally written with the two $\esup$ in \reff{optimale} switched. However, this can indeed be done. First, thanks to the comparison theorem (which still holds true for Lipschitz 2BSDEs, see \cite{stz}), we can conclude that $Y^{\ast}_t\geq Y^{u}_t$. Then since $u^{\ast}\in\mathcal{U}$, we have $Y^{\ast}_t\leq \underset{u\in\mathcal U}{\esup^{\mathbb P}}\ Y^{u}_t$. Hence, we do have $$Y^{\ast}_t= \underset{u\in\mathcal U}{\esup^{\mathbb P}}\ Y^{u}_t,\text{ which implies }\exp(\theta Y^{\ast}_t) =\underset{u\in\mathcal U}{\esup^{\mathbb P}}\exp(\theta Y^{u}_t).$$ Then, using the same type of arguments as above (see also \cite{mpz} for related results), we can show that $u^*$ is then an optimal strategy independent of the probability measure considered.
\end{Remark}

\section{2BSDEs and monotone approximations}\label{sec.approx}

Since the quasi-sure 2BSDE theory is still not completely mature, we feel that it is important to convey the idea that while some techniques from the classical theory can work for 2BSDEs (as in the previous section), some of the most important ones, namely the monotone approximation techniques, may fail.
Hence, the first part of this section is devoted to explaining why the now classical exponential transformation initiated by Kobylanski may not work in the general case, and more generally highlights the main reasons preventing us from using monotone approximations. In the second part however, we prove an existence result using another type of approximation which has the very desirable property to be stationary.

\subsection{Why the exponential transformation fails in general ?}
Coming back to Kobylanski \cite{kob}, we know that the exponential transformation used in the previous Section is an important tool in the study of quadratic BSDEs. However, unlike with a purely quadratic generator, in the general case the exponential change does not lead immediately to a Lipschitz BSDE. For the sake of clarity, let us consider the 2BSDE \reff{2bsdequad} and let us denote
$$\eta:=e^{\gamma\xi}, \text{ }\overline{Y}_t:=e^{\gamma Y_t}, \text{ }\overline{Z}_t:=\gamma \overline{Y}_tZ_t, \text{ }\overline{K}^{\mathbb P}_t:=\gamma\int_0^t\overline{Y}_sdK^{\mathbb P}_s-\sum_{0\leq s\leq t}e^{\gamma Y_s}-e^{\gamma Y_{s^-}}-\gamma\Delta Y_s e^{\gamma Y_{s^-}}.$$

Then we expect that, at least formally, if $(Y,Z)$ is a solution of \reff{2bsdequad}, then $(\overline{Y},\overline{Z})$ is a solution of the following $2$BSDE
\begin{equation}
\overline{Y}_t=\eta -\gamma\int_t^T\overline Y_s\left(\widehat F_s\left(\frac{\log \overline Y_s}{\gamma}, \frac{\overline Z_s}{\gamma\overline Y_s}\right)-\frac{\abs{\widehat a_s^{1/2}\overline Z_s}^2 }{2\gamma\overline Y_s^2}\right)ds-\int_t^T\overline Z_sdB_s+\overline K^\mathbb P_T-\overline K^\mathbb P_t.
\label{2bsdeapp}
\end{equation}

Let us now define for $y>0$, $G_t(\omega,y,z):=\gamma y\left(\widehat F_t\left(\omega,\frac{\log y}{\gamma}, \frac{z}{\gamma y}\right)-\frac{\abs{\widehat a_t^{1/2}z}^2 }{2\gamma y^2}\right).$ Then, despite the fact that the generator $G$ is not Lipschitz, it is possible, as shown by Kobylanski \cite{kob}, to find a sequence $(G^n)_{n\geq 0}$ of Lipschitz functions which decreases to $G$. Then, it is possible, thanks to the result of \cite{stz} to define for each $n$ the solution $(Y^n,Z^n)$ of the corresponding 2BSDE. The idea is then to prove existence and uniqueness of a solution for the 2BSDE with generator $G$ (and thus also for the 2BSDE \reff{2bsdequad}) by passing to the limit in some sense in the sequence $(Y^n,Z^n)$. Moreover, forgetting for a moment our particular quadratic framework, the exact same approach has been used countless times in the BSDE literature to tackle the wellposedness problem with different conditions on the generators.

\vspace{0.5em}
If we then follow the usual approach for standard BSDEs, the first step is to argue that thanks to the comparison theorem (which still holds true for Lipschitz 2BSDEs, see \cite{stz}), the sequence $Y^n$ is decreasing, and thanks to a priori estimates that it must converge $\mathcal P_H-q.s.$ to some process $Y$. And this is exactly now that the situation becomes almost inextricable with 2BSDEs. Indeed, if we were in the classical framework, this convergence of $Y^{n}$ together with the a priori estimates would be sufficient to prove the convergence in the usual $\mathbb H^2$ space, thanks to the dominated convergence theorem. However, in our case, since the norms involve the supremum over a family of probability measures which are mutually singular, this theorem generally fails (we refer the reader to Section $2.6$ in \cite{pos} for more details). More worrisome than that, it is generally believed that such a result is essentially wrong in any non-dominated setting. Therefore, we cannot be sure that
$$\underset{\mathbb P\in\mathcal P_H}{\sup}\mathbb E^\mathbb P\left[\int_0^T\abs{Y^n_t-Y_t}^2dt\right]\underset{n\rightarrow+\infty}{\longrightarrow}0.$$
Without this crucial step in the approximation proof, we are then stuck and cannot proceed any further.

\vspace{0.5em}
However, one could argue that only a monotone convergence Theorem is needed to obtain the desired result. To the best of our knowledge, the only monotone convergence Theorem in a similar setting has been proved by Denis, Hu and Peng (see \cite{dhp}). However, the set of probability measures has to be weakly compact, and one need to consider random variables $X^n$ which are regular in $\omega$, more precisely quasi-continuous, that is to say that for every $\eps>0$, there exists an open set $\mathcal O^\eps$ such that the $X^n$ are continuous in $\omega$ outside $\mathcal O^\eps$ and such that $$\underset{\mathbb P\in\mathcal P_H}{\sup}\mathbb P(\mathcal O^\eps)\leq\eps.$$ It is currently an open (an difficult) problem to know whether these conditions are sharp or not. Furthermore, it seems, at least to us, that it would be a very difficult task to weaken these assumptions to a point where the result could directly be applied to the sequence $(Y^n)_{n\geq 0}$ for the following reasons

\vspace{0.5em}
\hspace{0.9em}$\rm{(i)}$ \textit{First, if we assume that the terminal condition $\xi$ is in $UC_b(\Omega)$, since the generator $\widehat F$ (and thus $G^n$) are uniformly continuous in $\omega$, we can reasonably expect to be able to prove that the $Y^n$ will be also continuous in $\omega$, $\mathbb P-a.s.$, for every $\mathbb P\in\mathcal P_H$. However, this is clearly not sufficient to obtain the quasi-continuity. Indeed, for each $\mathbb P$, we would have a $\mathbb P$-negligible set outside of which the $Y^n$ are continuous in $\omega$. But since the probability measures are mutually singular, this does not imply the existence of the open set of the definition of quasi-continuity. We emphasize that this was to be expected because by definition, the solution of a 2BSDE is defined $\mathbb P-a.s.$ for every $\mathbb P$, and the quasi-continuity is by essence a notion related to the theory of capacities, not of probability measures.}

\vspace{0.5em}
\hspace{0.9em}$\rm{(ii)}$ \textit{Next, it has been shown that if we assume that the matrices $\underline a^\mathbb P$ and $\overline a^\mathbb P$ appearing in Definition \ref{defquad} are uniform in $\mathbb P$, then the set $\mathcal P_H$ is only weakly relatively compact. Then, we are left with two options. First, we can restrict the problem to a closed subset of $\mathcal P_H$, which will therefore be weakly compact. However, as pointed out in \cite{stz2}, it is not possible to restrict arbitrarily the probability measures considered. Indeed, since the whole approach of \cite{stz} to prove existence of Lipschitz 2BSDEs relies on stochastic control  and the dynamic programming equation, we need the set of processes $\alpha$ in the definition of $\overline{\mathcal P}_S$ (that is to say our set of control processes) to be stable by concatenation and bifurcation (see for instance Remark $3.1$ in \cite{bt2}) in order to recover the results of \cite{stz}. And it is not clear at all to us whether it is possible to find a closed subset of $\mathcal P_H$ satisfying this stability properties.}

\vspace{0.5em}
\textit{Otherwise, we could work with the weak closure of $\mathcal P_H$. This is exactly what is usually done in the so-called $G$-expectation literature initiated by Peng \cite{peng}. The problem now is that the probability measures in that closure no longer satisfy necessarily the martingale representation property and the 0-1 Blumenthal law. In that case (since the filtration $\mathbb F$ will only be quasi-left continuous), and as already shown by El Karoui and Huang \cite{ekh}, we would need to redefine a solution of a 2BSDE by adding a martingale orthogonal to the canonical process. However, defining such solutions is a complicated problem outside of the scope of this paper.}

\vspace{0.8em}
We hope to have convinced the reader that because of all the reasons listed above, monotone approximation should not be considered as the first method of choice in order to prove wellposedness results for 2BSDEs, and alternatives must be found. This is exactly what we do in Section \ref{section.1quad}, following the original approach of \cite{stz}. Notice nonetheless, that this by no way means that monotone approximations {\emph never} work. Indeed, in \cite{pos}, the author uses such an approach to prove existence of a solution to a 2BSDE with a generator with linear growth satisfying some monotonicity condition. The idea is that in this case it is possible to show that the sequence of approximated generators converges uniformly in $(y,z)$, and this allows to have a control on the difference $\abs{Y^n_t-Y_t}$ by a quantity which is regular enough to apply the monotone convergence Theorem of \cite{dhp}. Nonetheless, this relies heavily on the type of approximation used and cannot a priori be extended to more general cases.

\vspace{0.2em}
Notwithstanding this, we will show an existence result in the next subsection using an approximation which has the particularity of being stationary, which immediately solves the convergence problems that we mentioned above. This approach is based on very recent results of Briand and Elie \cite{be} on standard quadratic BSDEs.

\subsection{A stationary approximation}
\label{stationary_approx}

For technical reasons that we will explain below, we will work throughout this subsection under a subset of $\mathcal P_H$, which was first introduced in \cite{stz3}. Namely, we will denote by $\Xi$ the set of processes $\alpha$ satisfying
\begin{align}
 \nonumber&\alpha_t(\omega) = \sum_{n=0}^{+\infty} \sum_{i=1}^{+\infty} \alpha^{n,i}_t \mathbf{1}_{E^i_n}(\omega) \mathbf{1}_{[\tau_n(\omega),\tau_{n+1}(\omega))}(t),
 \label{doublesum}
\end{align}
where for each $i$ and for each $n$, $\alpha^{n,i}$ is a bounded deterministic mapping, $\tau_n$ is an $\F$-stopping time with $\tau_0=0$, such that $\tau_n < \tau_{n+1}$ on $\{\tau_n < +\infty \}$,  $\inf \{n \geq 0, \, \tau_n = +\infty\} <+\infty$, $\tau_n$ takes countably many values in some fixed $I_0\subset [0,T]$ which is countable and dense in $[0,T]$ and for each $n$, $(E_i^n)_{i \geq1} \subset \Fc_{\tau_n}$ forms a partition of $\Omega$.

\vspace{0.2em}
We will then consider the set $\widehat{\mathcal P}_H:=\left\{\mathbb P^\alpha\in\mathcal P_H,\ \alpha\in\Xi\right\}.$ As shown in \cite{stz2}, this set satisfies the right stability properties (already mentioned in the previous subsection) so much so that the Lipschitz theory of 2BSDEs still holds when we are working $\widehat{\mathcal P}_H-q.s.$ Notice that for the sake of simplicity, we will keep the same notations for the spaces considered under $\widehat{\mathcal P}_H$ or $\mathcal P_H$. Let us now describe the Assumptions under which we will be working
\begin{Assumption}\label{assump.hhh}
Let Assumption \ref{assump.hquad} holds, with the addition that the process $\phi$ in $\rm{(v)}$ is bounded and that the mapping $F$ is deterministic.
\end{Assumption}

The main result of this Section is then
\begin{Theorem}
\label{malliavin}
Let Assumption \ref{assump.hhh} hold. Assume further that $\xi\in\mathcal L^\infty_H$, that it is Malliavin differentiable $\widehat{\mathcal P}_H-q.s.$ and that its Malliavin derivative is in $\mathbb D^\infty_H$. Then the 2BSDE \reff{2bsdequad} (considered $\widehat{\mathcal P}_H-q.s.$) has a unique solution $(Y,Z)\in\mathbb D^\infty_H\times\mathbb H^2_H$. Moreover, the family $\{K^\mathbb P,\ \mathbb P\in\widehat{\mathcal P}_H\}$ can be aggregated.
\end{Theorem}

\proof
Uniqueness follows from Theorem \ref{quad.unique}, so we concentrate on the existence part. Let us define the following sequence of generators
$$F^n_t(y,z,a):=F_t\left(y,\frac{\abs{z}\wedge n}{\abs{z}}z,a\right),\text{ and }\widehat F^n_t(y,z):=F^n_t(y,z,\widehat a_t).$$

Then for each $n$, $F^n$ is uniformly Lipschitz in $(y,z)$ and thanks to Assumption \ref{assump.hhh}, we can apply the result of \cite{stz} to obtain the existence of a solution $(Y^n,Z^n)$ to the 2BSDE
\begin{equation}
Y^n_t=\xi+\int_t^T\widehat F^n_s(Y_s^n,Z_s^n)ds-\int_t^TZ_s^ndB_s+K_T^{\mathbb P,n}-K_t^{\mathbb P,n},\ \mathbb P-a.s.,\text{ for all }\mathbb P\in\widehat{\mathcal P}_H.
\label{eq:2BSDEn}
\end{equation}

Moreover, we have for all $\mathbb P\in\widehat{\mathcal P}_H$ and for all $t\in[0,T]$, $Y_t^n=\underset{\mathbb P^{'}\in\widehat{\mathcal P}_H(t^+,\mathbb P)}{\esup^\mathbb P}y_t^{\mathbb P,n}, \mathbb P-a.s.,$ where $(y^{\mathbb P,n},z^{\mathbb P,n})$ is the unique solution of the Lipschitz BSDE with generator $\widehat F^n$ and terminal condition $\xi$ under $\mathbb P$. Now, using Lemma $2.1$ in \cite{be} and its proof (see Remark \ref{rem.explain} below) under each $\mathbb P\in\widehat{\mathcal P}_H$, we know that the sequence $y^{\mathbb P,n}$ is actually stationary. Therefore, this also implies that the sequence $Y^n$ is stationary. Hence, we immediately have that $Y^n$ converges to some $Y$ in $\mathbb D^\infty_H$. Moreover, we still have the representation $Y_t=\underset{\mathbb P^{'}\in\widehat{\mathcal P}_H(t^+,\mathbb P)}{\esup^\mathbb P}y_t^{\mathbb P}, \mathbb P-a.s.$

\vspace{0.2em}
Now, identifying the martingale parts in \reff{eq:2BSDEn}, we also obtain that the sequence $Z^n$ is stationary and thus converges trivially in $\mathbb H^2_H$ to some $Z$. For $n$ large enough, we thus have $\widehat F^n_t(Y^n_t,Z^n_t)=\widehat F^n_t(Y_t,Z_t).$ Besides, we have by Assumption \ref{assump.hhh}
\begin{align*}
\abs{\widehat F^n_t(Y_t,Z_t)}&\leq \alpha +\beta \abs{Y_t}+\frac\gamma2\abs{\widehat a^{1/2}\frac{\abs{Z_t}\wedge n}{\abs{Z_t}}Z_t}^2\leq \alpha +\beta \abs{Y_t}+\frac\gamma2\abs{\widehat a^{1/2}Z_t}^2,\ \widehat{\mathcal P}_H-q.s.
\end{align*}

Since $(Y,Z)\in\mathbb D^\infty_H\times\mathbb H^2_H$, we can apply the dominated convergence theorem for the Lebesgue measure to obtain by continuity of $F$ that
\begin{equation*}
\int_0^T\widehat F^n_s(Y_s^n,Z_s^n)ds\underset{n\rightarrow+\infty}{\longrightarrow}\int_0^T\widehat F_s(Y_s,Z_s)ds, \widehat{\mathcal P}_H-q.s.
\end{equation*}

Using this result in \reff{eq:2BSDEn}, this implies necessarily that for each $\mathbb P$, $K^{\mathbb P,n}$ converges $\mathbb P-a.s.$ to a non-decreasing process $K^\mathbb P$. We can then check that the processes $K^\mathbb P$ satisfy the minimum condition \reff{2bsde.minquad} by proceeding exactly as in the proof of Proposition \ref{prop.pure}. Finally, the fact that the processes $K^\mathbb P$ can be aggregated is a direct consequence of the general aggregation result of Theorem $5.1$ in \cite{stz3}. Indeed, it is direct consequence of \reff{2bsde.minquad}, that if $\mathbb P^{'}\in\widehat{\mathcal P}_H(t^+,\mathbb P)$ for some $t$ and some $\mathbb P$, then we have $K_s^\mathbb P=K_s^{\mathbb P^{'}},\ 0\leq s\leq t,\ \mathbb P-a.s.$
\ep

\begin{Remark}\label{rem.explain}
We emphasize that the result of Lemma $2.1$ in \cite{be} can only be applied when the generator is deterministic. However, even though $F$ is indeed deterministic, $\widehat F$ is not, because $\widehat a$ is random. Nonetheless, given the particular form for the density of the quadratic variation of the canonical process we assumed in the definition of $\widehat{\mathcal P}_H$, we can apply the result of Briand and Elie between the stopping times and on each set of the partition of $\Omega$, since then $\widehat a$ and thus $\widehat F$ is indeed deterministic.
\end{Remark}

\section{A pathwise proof of existence}\label{section.1quad}

We have seen in the previous Section that it is usually inadequate to try and prove existence of a solution to a 2BSDE using monotone approximation techniques. Nonetheless, we have shown in Theorem \ref{quad.unique} that if a solution exists, it will necessarily verify the representation \reff{2bsde.minquad}. This gives us a natural candidate for the solution as a supremum of solutions to standard BSDEs. However, since those BSDEs are all defined on the support of mutually singular probability measures, it is not straightforward to define such a supremum, because of the problems raised by the negligible sets. In order to overcome this, Soner, Touzi and Zhang proposed in \cite{stz} a pathwise construction of the solution to a 2BSDE. Let us describe briefly their strategy.

\vspace{0.2em}
The first step is to define pathwise the solution to a standard BSDE. For simplicity, let us consider first a BSDE with a generator equal to $0$. Then, we know that the solution is given by the conditional expectation of the terminal condition. In order to define this solution pathwise, we can use the so-called regular conditional probability distribution (r.p.c.d. for short) of Stroock and Varadhan \cite{str}. In the general case, the idea is similar and consists on defining BSDEs on a shifted canonical space.

\vspace{0.2em}
Finally, we have to prove measurability and regularity of the candidate solution thus obtained, and the decomposition \reff{2bsdequad} is obtained through a non-linear Doob-Meyer decomposition. Our aim in this section is to extend this approach to the quadratic case. However, we want to insist on the fact that under our Assumption \ref{assump.hquad}, the proof is actually very close to the original one in \cite{stz}, which implies that we will only sketch it. Our main objective is this section is therefore to point out to the reader what are the important properties that the BSDEs considered must satisfy in order to be able to apply the proof strategy of \cite{stz}. Our contribution also resides on the technical result that we prove in Proposition \ref{prop.tech}, which extends the natural property of r.c.p.d. to the more general notion of BSDEs. We will see that our proof works in almost all cases of interest.

\subsection{Notations}

For the convenience of the reader, we recall below some of the notations introduced in \cite{stz}.

\vspace{0.2em}
$\bullet$ For $0\leq t\leq T$, denote by $\Omega^t:=\left\{\omega\in C\left([t,T],\mathbb R^d\right),\text{ }w(t)=0\right\}$ the shifted canonical space, $B^t$ the shifted canonical process, $\mathbb P_0^t$ the shifted Wiener measure and $\mathbb F^t$ the filtration generated by $B^t$.

\vspace{0.2em}
$\bullet$ For $0\leq s\leq t\leq T$ and $\omega\in \Omega^s$, define the shifted path $\omega^t\in \Omega^t$, $\omega^t_r:=\omega_r-\omega_t,\text{ }\forall r\in [t,T].$

\vspace{0.2em}
$\bullet$ For $0\leq s\leq t\leq T$ and $\omega\in \Omega^s$, $\widetilde \omega\in\Omega^t$ define the concatenation path $\omega\otimes_t\widetilde \omega\in\Omega^s$ by
$$(\omega\otimes_t\widetilde \omega)(r):=\omega_r1_{[s,t)}(r)+(\omega_t+\widetilde\omega_r)1_{[t,T]}(r),\text{ }\forall r\in[s,T].$$

\vspace{0.2em}
$\bullet$ For $0\leq s\leq t\leq T$ and a $\mathcal F^s_T$-measurable random variable $\xi$ on $\Omega^s$, for each $\omega \in\Omega^s$, define the shifted $\mathcal F^t_T$-measurable random variable $\xi^{t,\omega}$ on $\Omega^t$ by $\xi^{t,\omega}(\widetilde\omega):=\xi(\omega\otimes_t\widetilde \omega),\text{ }\forall \widetilde\omega\in\Omega^t.$ Similarly, for an $\mathbb F^s$-progressively measurable process $X$ on $[s,T]$ and $(t,\omega)\in[s,T]\times\Omega^s$, the shifted process $\left\{X_r^{t,\omega},r\in[t,T]\right\}$ is $\mathbb F^t$-progressively measurable.

\vspace{0.2em}
$\bullet$ For a $\mathbb F$-stopping time $\tau$, the r.c.p.d. of $\mathbb P$ (denoted $\mathbb P^\omega_\tau$) is a probability measure on $\mathcal F_T$ such that
$$\mathbb E_\tau^{\mathbb P}[\xi](\omega)=\mathbb E^{\mathbb P^\omega_\tau}[\xi],\text{ for }\mathbb P-a.e.\ \omega.$$

It also induces naturally a probability measure $\mathbb P^{\tau,\omega}$ (that we also call the r.c.p.d. of $\mathbb P$) on $\mathcal F_T^{\tau(\omega)}$ which in particular satisfies that for every bounded and $\mathcal F_T$-measurable random variable $\xi$
$$\mathbb E^{\mathbb P^\omega_\tau}\left[\xi\right]= \mathbb E^{\mathbb P^{\tau,\omega}}\left[\xi^{\tau,\omega}\right].$$

\vspace{0.2em}
$\bullet$ We define similarly as in Section \ref{section.1quad} the set $\overline{\mathcal P}^{t}_S$, by restricting to the shifted canonical space $\Omega^t$, and its subset $\mathcal P^{t}_H$.

\vspace{0.2em}
$\bullet$ Finally, we define our "shifted" generator, $\widehat F^{t,\omega}_s(\widetilde\omega,y,z):=F_s(\omega\otimes_t\widetilde\omega,y,z,\widehat a^t_s(\widetilde\omega)), \text{ }\forall (s,\widetilde\omega)\in[t,T]\times\Omega^t.$ Notice that thanks to Lemma $4.1$ in \cite{stz2}, this generator coincides for $\mathbb P$-a.e. $\omega$ with the shifted generator as defined above, that is to say $F_s(\omega\otimes_t\widetilde\omega,y,z,\widehat a_s(\omega\otimes_t\widetilde\omega)).$ The advantage of the chosen "shifted" generator is that it inherits the uniform continuity in $\omega$ under the $\mathbb L^\infty$ norm of $F$.

\subsection{Existence when $\xi$ is in $\rm{UC_b}(\Omega)$ }\label{Ucb}
To prove existence, as in \cite{stz}, we define the following value process $V_t$ pathwise
\begin{equation}
V_t(\omega):=\underset{\mathbb P\in\mathcal P^{t}_H}{\sup}\mathcal Y^{\mathbb P,t,\omega}_t\left(T,\xi\right), \text{ for all } \left(t,\omega\right)\in\left[0,T\right]\times\Omega,
\end{equation}

where, for any $\left(t_1,\omega\right)\in\left[0,T\right]\times\Omega,\ \mathbb P\in\mathcal P^{t_1}_H, \ t_2\in\left[t_1,T\right]$, and any $\mathcal F_{t_2}$-measurable $\eta\in\mathbb L^{\infty}\left(\mathbb P\right) $, we denote $\mathcal Y^{\mathbb P,t_1,\omega}_{t_1}\left(t_2,\eta\right):= y^{\mathbb P,t_1,\omega}_{t_1}$, where $\left(y^{\mathbb P,t_1,\omega},z^{\mathbb P,t_1,\omega}\right) $ is the solution of the following BSDE on the shifted space $\Omega^{t_1} $ under $\mathbb P$
\begin{equation}\label{eq.bsdeeee}
y^{\mathbb P,t_1,\omega}_{s}=\eta^{t_1,\omega}-\int^{t_2}_{s}\widehat{F}^{t_1,\omega}_{r}\left(y^{\mathbb P,t_1,\omega}_{r},z^{\mathbb P,t_1,\omega}_{r} \right)dr-\int^{t_2}_{s}z^{\mathbb P,t_1,\omega}_{r}dB^{t_1}_{r},\ s\in\left[t_1,t_2\right],\ \mathbb P-\text{a.s.}
\end{equation}

We recall that since the Blumenthal $0-1$ law holds for all our probability measures, $\mathcal Y^{\mathbb P,t,\omega}_{t}\left(1,\xi\right) $ is constant for any given $\left(t,\omega\right) $ and $\mathbb P\in\mathcal P^{t}_H $. Therefore, the process $V$ is well defined. However, we still do not know anything about its measurability. The following Lemma answers this question and explains the uniform continuity Assumptions in $\omega$.

\begin{Lemma}\label{unifcontquad}
Let $\xi\in \rm{UC_b}(\Omega)$. Under Assumption \ref{assump.hquad}, we have $\left|V_t\left(\omega\right)\right|\leq C\left(1+\No{\xi}_{\mathbb L^\infty_H}\right),$ $\text{for all }\left(t,\omega\right)\in\left[0,T\right]\times\Omega.$ Furthermore, $\left|V_t\left(\omega\right)-V_t\left(\omega'\right)\right|\leq C\rho\left(\No{\omega-\omega'}_t\right)$. In particular, $V_t$ is $\mathcal F_t$-measurable for every $t\in\left[0,T\right]$.
\end{Lemma}

The proof uses (as in \cite{stz}) a priori estimates for our class of BSDEs, which can be obtained exactly as in our proof of Theorem \ref{quad.unique}. The only addition is that we have to use the uniform continuity in $\omega$ that we assumed for both $\xi$ and $F$. Therefore we omit it. Then, we show the same dynamic programming principle as Proposition $4.7$ in \cite{stz2}

\begin{Proposition}
Let $\xi\in \rm{UC_b}(\Omega)$. Under Assumption \ref{assump.hquad}, we have for all $0\leq t_1<t_2\leq T$ and for all $\omega \in \Omega$
$$V_{t_1}(\omega)=\underset{\mathbb P\in \mathcal P^{t_1}_H}{\sup}\mathcal Y_{t_1}^{\mathbb P,t_1,\omega}(t_2,V_{t_2}^{t_1,\omega}).$$
\end{Proposition}

Once again, we can follow the proof in \cite{stz}, since it uses only comparison and stability results, which hold for our class of quadratic BSDEs. There is however one important technical point. Indeed, in their proof, they made use of the fact that solutions to Lipschitz BSDEs could be constructed via Picard iterations to justify that
\begin{equation}\label{rcpd}
y_t^\mathbb P(\omega)=\mathcal Y_t^{\mathbb P^{t,\omega},t,\omega}(T,\xi), \text{ for } \mathbb P-a.e.\text{ } \omega\in\Omega.
\end{equation}

In fact, since at each step of the iteration, the solution can be formulated as a conditional expectation under $\mathbb P$, this is a direct consequence of the properties of the r.c.p.d. However, as soon as we go away from the Lipschitz case, it is not clear at all whether solutions to BSDEs can be constructed this way (see nonetheless Tevzadze \cite{tev} where he uses a fixed-point argument to construct solutions to quadratic BSDEs but under stronger assumptions than ours), which implies that their argument no longer works directly in our case. The following Proposition aims at filling this gap, and shows that \reff{rcpd} holds whenever the corresponding BSDEs are wellposed\footnote{During the revision of this paper, Lin \cite{lin} proved a similar result in the quadratic case. However, his proof crucially uses the quadratic structure and cannot be extended to other cases. Our result is thus more general.}.

\begin{Proposition}\label{prop.tech}
Fix some $\mathbb P\in\mathcal P_H$ and some $t\in[0,T]$. Assume that the BSDEs \reff{bsdep} and \reff{eq.bsdeeee} are well posed. Then \reff{rcpd} holds.
\end{Proposition}

\proof
We start by proving a useful result. Let $\mathcal N^\mathbb P$ be a $\mathbb P$-null set and define for any $\omega\in\Omega$, the following subset of $\Omega^t$, $\mathcal N^{\mathbb P,t,\omega}:=\left\{\widetilde\omega\in\Omega^t, \ \omega\otimes_t\widetilde\omega\in\mathcal N^\mathbb P\right\}.$ We claim that for $\mathbb P$-a.e. $\omega$, $\mathcal N^{\mathbb P,t,\omega}$ is a $\mathbb P^{t,\omega}$-null set. Indeed, we have by definition
$$1_{\mathcal N^{\mathbb P,t,\omega}}(\widetilde\omega)=1_{\mathcal N^\mathbb P}(\omega\otimes_t\widetilde\omega)=\left(1_{\mathcal N^\mathbb P}\right)^{t,\omega}(\widetilde\omega).$$

Hence, by definition of the r.c.p.d., we have for $\mathbb P$-a.e. $\omega$
\begin{align*}
\mathbb P^{t,\omega}\left(\mathcal N^{\mathbb P,t,\omega}\right)=\mathbb E^{\mathbb P^{t,\omega}}\left[\left(1_{\mathcal N^\mathbb P}\right)^{t,\omega}\right]=\mathbb E^\mathbb P_t\left[1_{\mathcal N^\mathbb P}\right](\omega),
\end{align*}
which is equal to $0$ for $\mathbb P$-a.e. $\omega$ since $\mathcal N^\mathbb P$ is a $\mathbb P$-null set. Hence the desired result. The proof will now be done in two steps.

\vspace{0.2em}
\textbf{Step (i)}: Let us denote $\mathcal N^{\mathbb P}$ the $\mathbb P$-null set outside of which (\ref{bsdep}) holds. We emphasize that this set will be enlarged later on with other $\mathbb P$-null sets, but we will keep the same notation to denote it for the sake of simplicity. Then, for any $t\in\left[0,T\right]$ and $\omega\in\Omega/\mathcal N^{\mathbb P}$, we define $\mathcal N^{\mathbb P,\omega,t}$ as above, which is a $\mathbb P^{t,\omega}$-null set for $\mathbb P$-a.e. $\omega$. Then for $\tilde{\omega}\in\Omega^t/\mathcal N^{\mathbb P,\omega,t}$, we have
\begin{align}\label{eq.bsdeeee2}
\nonumber y^{\mathbb P}_{s}(\omega\otimes_t\widetilde{\omega})=&\ \xi(\omega\otimes_t\widetilde{\omega})-\int^{T}_{s}F_{r}\left(\omega\otimes_t\widetilde{\omega},y^{\mathbb P}_{r}(\omega\otimes_t\widetilde{\omega}),z^{\mathbb P}_{r}(\omega\otimes_t\widetilde{\omega}),\widehat{a}_r(\omega\otimes_t\widetilde{\omega}) \right)dr\\
&-\left(\int^{T}_{s}z^{\mathbb P}_{r}dB_{r}\right)(\omega\otimes_t\widetilde{\omega}),\ s\in\left[t,T\right].
\end{align}

Remember now that by Lemma $4.1$ of \cite{stz2} we have for $\mathbb P$-a.e. $\omega$, $\widehat a_r(\omega\otimes_t\widetilde{\omega})=\widehat a^t_r(\widetilde\omega),\ \widetilde\omega\in\Omega^t\backslash\mathcal N^{\mathbb P^{t,\omega}},$ where $\mathcal N^{\mathbb P^{t,\omega}}$ is another $\mathbb P^{t,\omega}$-null set.
Thus, enlarging the $\mathbb P$-null set $\mathcal N^\mathbb P$ to a larger $\mathbb P$-null set if necessary, for any $(\omega,\widetilde\omega)\in\left(\Omega\backslash\mathcal N^\mathbb P\right)\times\left(\Omega^t\backslash\left\{\mathcal N^{\mathbb P,\omega,t}\cup\mathcal N^{\mathbb P^{t,\omega}}\right\}\right)$, we have
\begin{align*}
y^{\mathbb P}_{s}(\omega\otimes_t\widetilde{\omega})=&\ \xi^{\omega,t}(\widetilde{\omega})-\int^{T}_{s}\widehat F^{\omega,t}_{r}\left(\widetilde{\omega},y^{\mathbb P}_{r}(\omega\otimes_t\widetilde{\omega}),z^{\mathbb P}_{r}(\omega\otimes_t\widetilde{\omega}) \right)dr-\left(\int^{T}_{s}z^{\mathbb P}_{r}dB_{r}\right)(\omega\otimes_t\widetilde{\omega}).
\end{align*}

In Step (ii) below, we will prove that for $\mathbb P$-a.e. $\omega$ and for $\mathbb P^{t,\omega}$-a.e. $\widetilde\omega$, we have
\begin{align}\label{eq.sto}
\left(\int^T_sz^{\mathbb P}_{r}dB_{r}\right)(\omega\otimes_t\widetilde{\omega})=\left(\int^{T}_{s}\left(z^{\mathbb P}\right)^{t,\omega}_{r}dB^t_{r}\right)(\widetilde{\omega}).
\end{align}

Then the above equation means that $(y^{\mathbb P}_{s}(\omega\otimes_t\cdot), z^{\mathbb P}_{s}(\omega\otimes_t\cdot))$ is a solution of \reff{eq.bsdeeee} under $\mathbb P^{t,\omega}$. Uniqueness of the solution, more precisely when taking $s=t$, gives the desired result.

\vspace{0.2em}
\textbf{Step (ii)}: We still need to prove that \reff{eq.sto} holds. Actually, the stochastic integral on the left-hand side of \reff{eq.sto} is defined outside some $\mathbb P$-null set denoted $\mathcal N^{\mathbb P}_1$, and the one on the right-hand side outside a $\mathbb P^{t,\omega}$-null set denoted $\mathcal N^{\mathbb P^{t,\omega}}_1$. Then, defining as in the previous step the set $\mathcal N^{\mathbb P,t,\omega}_1$, we know that the left-hand side is well-defined for any $\omega$ outside $\mathcal N^\mathbb P_1$ and any $\widetilde\omega$ outside $\mathcal N^{\mathbb P,t,\omega}_1$. Moreover, by definition of the stochastic integral, we can find a well-chosen sequence $(\{r_i^n\}_{0\leq i\leq n})_{n\geq 0}$ of partitions of $[s,T]$ and such that
$$\left(\int^T_sz^{\mathbb P}_{r}dB_{r}\right)(\omega)=\underset{n\rightarrow +\infty}{\lim}\ I^n(\omega),\text{ for $\omega\in\Omega\backslash\mathcal N^\mathbb P_1$,}$$
where
$$I^n(\omega):=\sum_{i=0}^nz^{n,\mathbb P}_{r_i^n}(\omega)\left(B_{r_{i+1}^n}(\omega)-B_{r^n_i}(\omega)\right),$$
where $z^{n,\P}$ can be any sequence of simple processes converging to $z^\P$ in $\H^2(\P)$. Then, we have by definition of the canonical process
$$I^n(\omega\otimes_t\widetilde\omega)=\sum_{i=0}^n\left(z^{n,\mathbb P}\right)^{t,\omega}_{r_i^n}(\widetilde\omega)\left(B^t_{r_{i+1}^n}(\widetilde\omega)-B^t_{r^n_i}(\widetilde\omega)\right),$$
which converges in probability for $\mathbb P^{t,\omega}$ to the right-hand side \reff{eq.sto} (up to a subsequence if necessary). Indeed, it is not difficult to see that if $z^{n,\P}$ converges to $z^\P$ in $\H^2(\P)$, then $(z^{n,\P})^{t,\omega}$ converges to $(z^\P)^{t,\omega}$ in $\H^2(\P^{t,\omega})$ for $\P$-a.e. $\omega$, along a subsequence if necessary. Choosing another subsequence, this convergence holds $\mathbb P^{t,\omega}$-a.s., and the previous one still holds $\mathbb P$-a.s. along this subsequence. Finally, since all the sets considered here are either $\mathbb P$-null sets or $\mathbb P^{t,\omega}$-null sets (for $\mathbb P$-a.e. $\omega$), the result follows.
\ep

\vspace{0.2em}
Now that we solved the measurability issues for $V_t$, we need to study its regularity in time. However, it seems difficult to obtain a result directly, given the definition of $V$. This is the reason why we define now for all $(t,\omega)$, the $\mathbb F^+$-progressively measurable process
$$V_t^+:=\underset{r\in\mathbb Q\cap(t,T],r\downarrow t}{\overline \lim}V_r.$$

This new value process will then be proved to be c\`adl\`ag. Notice that a priori $V^+$ is only $\mathbb F^+$-progressively measurable, and not $\mathbb F$-progressively measurable. This explains why in the definition of the spaces in Section \ref{sec.space}, the processes are assumed to be $\mathbb F^+$-progressively measurable.

\begin{Lemma}
Under the conditions of the previous Proposition, we have
$$V_t^+=\underset{r\in\mathbb Q\cap(t,T],r\downarrow t}{\lim}V_r,\text{ }\mathcal P_H-q.s.$$
and thus $V^+$ is c\`adl\`ag $\mathcal P_H-q.s.$
\end{Lemma}

\proof
Actually, we can proceed exactly as in the proof of Lemma $4.8$ in \cite{stz2}, since the theory of $g$-expectations of Peng has been extended by Ma and Yao in \cite{ma} to the quadratic case (see in particular their Corollary $5.6$ for our purpose).
\ep

\vspace{0.2em}
Finally, proceeding as in Steps $1$ and $2$ of the proof of Theorem $4.5$ in \cite{stz2}, and in particular using the Doob-Meyer decomposition proved in \cite{ma} (Theorem $5.2$), we can get the existence of a process $Z$ and a family of non-decreasing processes $\left\{K^\mathbb P,\mathbb P\in\mathcal P_H\right\}$ such that
$$V_t^+=V_0^++\int_0^t\widehat F_s(V_s^+,Z_s)ds+\int_0^tZ_sdB_s-K_t^\mathbb P, \text{ }\mathbb P-a.s. \text{ }\forall \mathbb P\in\mathcal P_H.$$

For the sake of completeness, we provide the representation \reff{representationquad} for $V$ and $V^+$, and that, as shown in Proposition $4.11$ of \cite{stz2}, we actually have $V=V^+$, $\mathcal P_H-q.s.$, which shows that in the case of a terminal condition in $UC_b(\Omega)$, the solution of the $2$BSDE is actually $\mathbb F$-progressively measurable. This will be important in Section \ref{section.3quad}.

\begin{Proposition}\label{prop.rep}
Let $\xi\in UC_b(\Omega)$. Under Assumption \ref{assump.hquad}, $V_t=V_t^+, \text{ }\mathcal P_H-q.s.$ and
$$V_t=\underset{\mathbb P^{'}\in\mathcal P_H(t,\mathbb P)}{\esup^\mathbb P}\mathcal Y_t^{\mathbb P^{'}}(T,\xi)\text{ and } V_t^+=\underset{\mathbb P^{'}\in\mathcal P_H(t^+,\mathbb P)}{\esup^\mathbb P}\mathcal Y_t^{\mathbb P^{'}}(T,\xi), \text{ }\mathbb P-a.s., \text{ }\forall \mathbb P\in\mathcal P_H.$$
\end{Proposition}

To be sure that we have found a solution to our $2$BSDE, it remains to check that the family of non-decreasing processes above satisfies the minimum condition. However, this can be done exactly as in the proof of Proposition \ref{prop.pure}.

\subsection{Main result}
We are now in position to state the main result of this section which can be proved as in \cite{stz}, using the a priori estimates obtained in Theorem \ref{estimatesquad}. For the aggregation result, we only give of possible choice of assumptions and refer the reader to \cite{mpz2} for more details.

\begin{Theorem}\label{mainquad}
Let $\xi\in\mathcal L^{\infty}_H$. Under Assumption \ref{assump.hquad}, there exists a unique solution $(Y,Z)\in\mathbb D^{\infty}_H\times\mathbb H^{2}_H$ of the $2\rm{BSDE}$ \reff{2bsdequad}. If in addition we work in the ZFC model of set theory with the addition of the axiom of choice and the continuum hypothesis, then the corresponding family of processes $(K^\mathbb P)_{\mathbb P\in\mathcal P_H}$ can be aggregated into a universal process $K$.
\end{Theorem}

\begin{Remark}
As mentioned in the introduction, there is another proof of existence for classical quadratic BSDEs using a fixed point argument, which was obtained by Tevzadze \cite{tev}. However, this type of proof does not seem to work for the simple reason that it requires to have a general representation theorem for non-linear martingales (or $G$-martingales). But the only existing results where obtained in \cite{stz4} and \cite{psz} and require strong regularity assumptions in $\omega$ which are not verified by our generator, since it contains $\widehat a$ which has no regularity at all.
\end{Remark}

\section{Connection with fully non-linear PDEs}\label{section.3quad}
In this section, we assume that all the randomness in $H$ only depends on the current value of the canonical process $B$ (the so-called Markov property), $H_t(\omega,y,z,\gamma)=h(t,B_t(\omega),y,z,\gamma),$ where $h:[0,T]\times\mathbb R^d\times\mathbb R\times\mathbb R^d\times D_h\rightarrow\mathbb R$ is a deterministic map. Then, we define as in Section \ref{section.1quad} the corresponding conjugate and bi-conjugate functions
\begin{align*}
f(\cdot,a):=\underset{\gamma\in D_h}{\Sup}\left\{\frac12\Tr{a\gamma}-h(\cdot,\gamma)\right\},\ \widehat h(\cdot,\gamma):=\underset{a\in \mathbb S_d^{>0}}{\Sup}\left\{\frac12\Tr{a\gamma}-f(\cdot,a)\right\}.
\end{align*}

We denote $\mathcal P_h:=\mathcal P_H$, and following \cite{stz}, we strengthen Assumptions \ref{assump.hquad} 
\begin{Assumption}\label{assump.h4}
\begin{itemize}
\item[\rm{(i)}] The domain $D_{f_t}$ of $f$ in $a$ is independent of $(x,y,z)$.
\item[\rm{(ii)}] On $D_{f_t}$, f is uniformly continuous in $t$, uniformly in $a$.
\item[\rm{(iii)}] $f$ is continuous in $z$ and there exists $(\alpha, \beta,\gamma)$ such that
$$\abs{f(t,x,y,z,a)}\leq \alpha+\beta\abs{y}+\frac\gamma2|a^{1/2}z|^2,\text{ for all } (t,x,y,z,a).$$
\item[\rm{(iv)}] $f$ is uniformly continuous in $x$, uniformly in $(t,y,z,a)$, with a modulus of continuity $\rho$ which has polynomial growth.
\item[\rm{(v)}] There exists $\mu>0$ and a bounded $\mathbb R^d$-valued function $\phi$ such that for all $(t,y,z,z',a)$
$$|f(t,x,y,z,a)-f(t,x,y,z',a)-\phi(t).a^{1/2}(z-z')|\leq \mu a^{1/2}\abs{z-z'}\left(\abs{ a^{1/2}z}+\abs{ a^{1/2}z'}\right).$$
\item[\rm{(vi)}] $f$ is Lipschitz in $y$, uniformly in $(t,x,z,a)$.
\end{itemize}
\end{Assumption}

\vspace{0.2em}
Let now $g:\mathbb R^d\rightarrow \mathbb R$ be a Lebesgue measurable and bounded function. Our object of interest here is the following $2$BSDE with terminal condition $\xi=g(B_T)$
\begin{equation}\label{2bsdemark}
Y_t=g(B_T)-\int_t^Tf(s,B_s,Y_s,Z_s,\widehat a_s)ds-\int_t^TZ_sdB_s+K^\mathbb P_T-K^\mathbb P_t,\text{ }\mathcal P_h-q.s.
\end{equation}

The aim of this section is to generalize the results of \cite{stz} and obtain the connection $Y_t=v(t,B_t)$, $\mathcal P_h-q.s.$, where $v$ verifies in some sense the following fully non-linear PDE
\begin{equation}\label{pde}
\frac{\partial v}{\partial t}(t,x)+\widehat h\left(t,x,v(t,x),Dv(t,x),D^2v(t,x)\right)=0,\text{ }t\in[0,T)\\[0.5em]
\end{equation}

Following the classical terminology in the BSDE literature, we say that the solution of the $2$BSDE is Markovian if it can be represented by a deterministic function of $t$ and $B_t$. In this subsection, we will construct such a function following the same spirit as in the construction in the previous section. Exactly as in the previous section, all the techniques used in \cite{stz} still work, because they only use comparison and stability results for BSDEs, which still hold in our case. We want to insist once more on this fact, which roughly says that as long as once considers BSDEs such that these properties are verified, then the pathwise construction of the 2BSDEs and their link with fully non-linear PDEs can be obtained exactly as in \cite{stz} (with of course some minor modification). For this reason, and given the length of the paper, we will therefore only state the results and refer the reader to \cite{stz} and the proofs of comparison theorems and a priori estimates in the present paper for more details.

\vspace{0.2em}
With the same notations for shifted spaces, we define for any $(t,x)\in[0,T]\times\mathbb R^d$, $B^{t,x}_s:=x+B^t_s,\ \text{for all }s\in[t,T].$ Let now $\tau$ be an $\mathbb F^t$-stopping time, $\mathbb P\in\mathcal P^t_h$ and $\eta$ a $\mathbb P$-bounded $\mathcal F_\tau^t$-measurable random variable. Similarly as in \reff{eq.bsdeeee}, we denote $(y^{\mathbb P,t,x}, z^{\mathbb P,t,x})$ the unique solution of the following BSDE
\begin{equation}\label{bsde.markov}
y_s^{\mathbb P,t,x}=\eta -\int_s^\tau f(u,B_u^{t,x},y_u^{\mathbb P,t,x},z_u^{\mathbb P,t,x},\widehat a^t_u)du-\int_s^\tau z_u^{\mathbb P,t,x}dB^{t,x}_u, \text{ }t\leq s\leq \tau,\text{ }\mathbb P-a.s.
\end{equation}

Next, we define the following deterministic function (by virtue of the Blumenthal $0-1$ law)
\begin{equation}\label{eq.valuef}
u(t,x):=\underset{\mathbb P\in\mathcal P^t_h}{\sup}\mathcal Y^{\mathbb P,t,x}_t(T,g(B^{t,x}_T)), \text{ for }(t,x)\in[0,T]\times\mathbb R^d.
\end{equation}

We then have the following Theorem, which is actually Theorem $5.9$ of \cite{stz} in our framework

\begin{Theorem}
Let Assumption \ref{assump.h4} hold, and assume that $g$ is bounded and uniformly continuous. Then the $2$BSDE \reff{2bsdemark} has a unique solution $(Y,Z)\in\mathbb D^{\infty}_H\times\mathbb H^2_H$ and we have $Y_t=u(t,B_t)$. Moreover, $u$ is uniformly continuous in $x$ and right-continuous in $t$.
\end{Theorem}

\subsection{Non-linear Feynman-Kac formula in the quadratic case}

Exactly as in the classical case and as in Theorem $5.3$ in \cite{stz}, we have a non-linear version of the Feynman-Kac formula. The proof is the same as in \cite{stz}, so we omit it. Notice however that it is more involved than in the classical case, mainly due to the technicalities introduced by the quasi-sure framework.

\begin{Theorem}
Under Assumption \ref{assump.h4}, suppose that $\widehat h$ is continuous in its domain, that $D_f$ is independent of $t$ and is bounded both from above and away from $0$. Let $v\in C^{1,2}([0,T),\mathbb R^d)$ be a classical solution of \reff{pde} with $\left\{(v,Dv)(t,B_t)\right\}_{0\leq t\leq T}\in\mathbb D^{\infty}_H\times\mathbb H^2_H$. Then
$$Y_t:=v(t,B_t),\text{ }Z_t:=Dv(t,B_t),\text{ }dK_t:=k_sds,$$
is the unique solution of the quadratic $2$BSDE \reff{2bsdemark}, where
$$k_t:=\widehat h(t,B_t,Y_t,Z_t,\Gamma_t)-\frac12\Tr{\widehat a_t^{1/2}\Gamma_t}+f(t,B_t,Y_t,Z_t,\widehat a_t)\text{ and }\Gamma_t:=D^2v(t,B_t).$$
\end{Theorem}

\subsection{The viscosity solution property}
As usual when dealing with possibly discontinuous viscosity solutions, we introduce the following upper and lower-semicontinuous envelopes
\begin{align*}
u_*(\vartheta):=\underset{\vartheta'\rightarrow\vartheta}{\underline{\lim}}u(\vartheta'),\text{ }u^*(\vartheta):=\underset{\vartheta'\rightarrow\vartheta}{\overline{\lim}}u(\vartheta'),\ \widehat h_*(\vartheta):=\underset{\scriptstyle\vartheta'\rightarrow\vartheta}{\underline{\lim}}\widehat h(\vartheta')\text{, }\widehat h^*(\vartheta):=\underset{\vartheta'\rightarrow\vartheta}{\overline{\lim}}\widehat h(\vartheta')
\end{align*}

In order to prove the main Theorem of this subsection, we will need the following Proposition, whose proof (which is rather technical) is omitted, since it is exactly the same as the proof of Propositions $5.10$ and $5.14$ and Lemma $6.2$ in \cite{stz}.

\begin{Proposition}\label{prop.dynamicp}
Let Assumption $\ref{assump.h4}$ hold. Then for any bounded function $g$ and $(t,x)$

\vspace{0.2em}
(i) $\forall$  $\left\{\tau^\mathbb P,\mathbb P\in\mathcal P_h^t\right\}$, $\mathbb F^t$-stopping times , we have $u(t,x)\leq\underset{\mathbb P\in \mathcal P^t_h}{\sup}\mathcal Y^{\mathbb P,t,x}_t(\tau^\mathbb P,u^*(\tau^\mathbb P,B^{t,x}_{\tau^\mathbb P})).$

\vspace{0.2em}
(ii) If in addition $g$ is lower-semicontinuous, then $u(t,x)=\underset{\mathbb P\in \mathcal P^t_h}{\sup}\mathcal Y^{\mathbb P,t,x}_t(\tau^\mathbb P,u(\tau^\mathbb P,B^{t,x}_{\tau^\mathbb P})).$
\end{Proposition}

Now we can state the main Theorem of this section
\begin{Theorem}
Let Assumption \ref{assump.h4} hold true. Then
\begin{itemize}
\item[\rm{(i)}] $u$ is a viscosity subsolution of $-\partial_tu^*-\widehat h^*(\cdot,u^*,Du^*,D^2u^*)\leq 0,\text{ on }[0,T)\times\mathbb R^d.$

\item[\rm{(ii)}] If in addition $g$ is lower-semicontinuous and $D_f$ is independent of $t$, then $u$ is a viscosity supersolution of $-\partial_tu_*-\widehat h_*(\cdot,u_*,Du_*,D^2u_*)\geq 0,\text{ on }[0,T)\times\mathbb R^d.$
\end{itemize}
\end{Theorem}

\begin{Remark}
In order to conclude one would of course need a comparison theorem to hold for viscosity solutions to the PDE \reff{pde}. In this respect, the quadratic growth assumption is not a major problem, and comparison will hold under additional assumptions on the functions $h$ and $f$. We refer the interested reader to Theorems $3.3$ and $5.1$ in \cite{crand} (and the references therein) for some examples of conditions.
\end{Remark}


\begin{thebibliography}{aa12}

 \bibitem{bar}
 Barrieu, P., El Karoui, N. (2008). Pricing, hedging and designing derivatives with risk measures, {\sl Chap. 3 in "Indifference Pricing: Theory and Applications", Springer}, 77--146.
 \bibitem{elkarbar} Barrieu, P., El Karoui, N. (2011). Monotone stability of quadratic
semimartingales with applications to general quadratic BSDEs, {\sl Ann. Prob.}, to appear.
\bibitem{bie}
Bielecki, T.R., Pliska, S.R. (1999). Risk-sensitive dynamic asset management, {\sl App. Math. Optim.}, 39:337--360.
\bibitem{bis} Bismut, J.M. (1973).
Conjugate convex functions in optimal stochastic control, {\sl J. Math. Anal. Appl.}, 44:384--404.
\bibitem{bt}
Bouchard, B., Touzi, N. (2004). Discrete-time approximation and Monte Carlo simulation of backward stochastic differential equations, {\sl Stoch. Proc. App.}, 111:175--206.
\bibitem{bt2}
Bouchard, B., Touzi, N. (2011). Weak Dynamic Programming Principle for Viscosity Solutions, {\sl SIAM J. Cont. Opt.}, 49(3):948--962.
\bibitem{bh} Briand, Ph., Hu, Y. (2006).
BSDE with quadratic growth and unbounded terminal value, {\sl Probab. Theory Relat. Fields}, 136:604--618.
\bibitem{bh2} Briand, Ph., Hu, Y. (2008).
Quadratic BSDEs with convex generators and unbounded terminal conditions, {\sl Probab. Theory Relat. Fields}, 141:543--567.
\bibitem{be}
Briand, P., Elie, R. (2012). A simple constructive approach to quadratic BSDEs with or without delay, preprint.
\bibitem{cstv} Cheridito, P., Soner, H.M.,
Touzi, N., and Victoir, N. (2007). Second Order Backward Stochastic
Differential Equations and Fully Non-Linear Parabolic PDEs, {\sl
Comm. Pure App. Math.},  60(7):1081--1110.
\bibitem{crand}
Crandall, M.G., Ishii, H., Lions, P.-L. (1992). User's guide to viscosity solutions of second order PDEs, {\sl Bull. Am. Math. Soc.} 27(1):1--67.
\bibitem{denis} Denis, L., Martini, C. (2006). A theoretical framework for the pricing of contingent claims in the presence of model uncertainty, {\sl Ann. of App. Prob.}, 16(2):827--852.
\bibitem{dhp} Denis, L., Hu, M., Peng, S. (2011). Function spaces and capacity related to a
sublinear expectation: application to
G-Brownian motion paths, {\sl Pot. Anal.}, 34(2):139--161.
\bibitem{elkaroui}
El Karoui, N., Peng, S. and Quenez, M.C. (1994). Backward stochastic differential equations in finance, {\sl Mathematical Finance}, 7(1):1--71.
\bibitem{ekh}
El Karoui, N., Huang, S.-J. (1997). A general result of existence and uniqueness of backward stochastic differential equation, {\sl in Backward Stochastic Differential Equations, Pitman Research Notes}, 364:27--36.
\bibitem{elkarmat}
El Karoui, N., Hamad\`ene, S., Matoussi, A. (2008). Backward stochastic differential equations and applications, {\sl Chapter 8 in "Indifference Pricing: Theory and Applications", Springer-Verlag}, 267--320.
\bibitem{flem00}
Fleming, W.H., Sheu, S.J. (2000). Risk-sensitive control and optimal investment model, {\sl Math. Fin.}, 10(17):197--213.
\bibitem{flem02}
Fleming, W.H., Sheu, S.J. (2002). Risk-sensitive control and optimal investment model II, {\sl Ann. App. Prob.}, 12(2):730--767.
\bibitem{kar} Karandikar, R. (1995). On pathwise stochastic integration, {\sl Stoc. Proc. Ap.}, 57:11--18.
\bibitem{kaz} Kazamaki, N. (1994). Continuous exponential martingales and BMO. {\sl Springer}.
\bibitem{kob}
Kobylanski, M. (2000). Backward stochastic differential equations and partial differential equations with quadratic growth, {\sl Ann. Prob.} 28:259--276.
\bibitem{lin}
Lin, Y. (2013). A new result for second order BSDEs with quadratic growth and its applications, preprint, {\sl  arXiv:1301.0457}.
\bibitem{ma}
Ma, J., Yao, S. (2010). On quadratic g-evaluations$\slash$ expectations and related analysis, {\sl Stoch. Ana. and App.}, 28(4):711-734.
\bibitem{mpz}
Matoussi, A., Possama\"{i}, D., Zhou, C. (2011). Robust utility maximization under volatility uncertainty, {\sl Math. Fin.}, to appear.
\bibitem{mpz2}
Matoussi, A., Possama\"{i}, D., Zhou, C. (2011). Second-order reflected BSDEs, {\sl Ann. App. Prob.}, to appear.
\bibitem{nutz}
Nutz, M. (2012). Pathwise construction of stochastic integrals, {\sl Elec. Com. Prob.}, 17(24):1--7.
\bibitem{pardpeng}
Pardoux, E. and Peng, S (1990). Adapted solution of a backward stochastic differential equation, {\sl Systems Control Lett.}, 14:55--61.
\bibitem{pardpeng2}
Pardoux, E., Peng, S (1992). Backward stochastic differential equations and quasilinear parabolic differential equations, {\sl Lecture notes in CIS}, 176:200--217.
\bibitem{peng} Peng, S. (2010). Nonlinear expectations and stochastic calculus under uncertainty, preprint.
\bibitem{psz}
Peng, S., Song, Y., Zhang, J. (2012). A complete representation theorem for $G$-martingales, preprint, {\sl  arXiv:1201.2629.}
\bibitem{pos} Possama\"{i}, D. (2010).
Second order backward stochastic differential equations with continuous coefficient, {\sl Stoc. Proc. App.}, to appear.

\bibitem{stz}
Soner, H.M., Touzi, N., Zhang J. (2012). Wellposedness of second order BSDE's, {\sl Prob. Th. and  Related Fields}, 153(1-2):149--190.
\bibitem{stz2}
Soner, H.M., Touzi, N., Zhang J. (2010). Dual formulation of second order target problems, {\sl Ann. of App. Prob.}, to appear.
\bibitem{stz4}
Soner, H.M., Touzi, N., Zhang J. (2010). Martingale representation theorem for the G-expectation, {\sl Stoch. Proc. and their App.}, 121:265--287.
\bibitem{stz3}
Soner, H.M., Touzi, N., Zhang J. (2011). Quasi-sure stochastic analysis through aggregation, {\sl Elec. Journal of Prob.}, 16(67):1844-1879.
\bibitem{str}
Stroock, D.W., Varadhan, S.R.S. (1979). Multidimensional diffusion processes, {\sl Springer.}
\bibitem{tev}
Tevzadze, R. (2008). Solvability of backward stochastic differential equations with quadratic growth, {\sl Stoch. Proc. and their App.}, 118:503--515.
\end{thebibliography}
\end{document}